\documentclass[letterpaper,10pt]{article}
\usepackage{amssymb}
\usepackage{graphicx}
\usepackage{epstopdf}
\usepackage{url,color}
\usepackage[algo2e,linesnumbered, vlined,ruled]{algorithm2e}
\usepackage{amsmath}
\usepackage{multirow}
\usepackage{float}

\def\d{~\mathrm{d}}
\def\x{\mathbf{x}}
\def\w{\mathbf{w}}
\def\j{\mathbf{j}}
\def\dx{~\mathrm{d}\mathbf{x}}
\def\b{\mathrm{b}}
\def\Tr{\mathrm{Tr}}

\newtheorem{algorithm}{{\bf Algorithm}}
\newtheorem{definition}{{\bf Definition}}
\newtheorem{conjecture}{{\bf Conjecture}}

\newtheorem{property}{{\bf Property}}

\begin{document}
\title{Compressed Modes for Variational Problems in Mathematics and Physics}
\date{}
\author{Vidvuds Ozoli\c{n}\v{s} \thanks{Department of Materials Science and Engineering, University of California, Los Angeles, USA.
                      ({\tt vidvuds@ucla.edu})}
           \and Rongjie Lai\thanks{Department of Mathematics, University of California, Irvine, USA.
                     ({\tt rongjiel@math.uci.edu})}
           \and Russel Caflisch\thanks{Department of Mathematics, University of California, Los Angeles, USA. 
                         ({\tt rcaflisch@ipam.ucla.edu}).} 
           \and Stanley Osher\thanks{Department of Mathematics, University of California, Los Angeles, USA. 
                         ({\tt sjo@math.ucla.edu}).} 
     }

\maketitle

\begin{abstract}
This paper describes a general formalism for obtaining localized solutions to a class of problems in mathematical physics, which can be recast as variational optimization problems. This class includes the important cases of Schr\"odinger's equation in quantum mechanics and electromagnetic equations for light propagation in photonic crystals.
These ideas can also be applied to develop a spatially localized basis that spans the eigenspace of a differential operator, for instance, the Laplace operator, generalizing the concept of plane waves to an orthogonal real-space basis with multi-resolution capabilities.
\end{abstract}

%\begin{keywords}Schr\"odinger's equation, Compressed Modes, Compressed Plane Waves.
%\end{keywords}

\section{Introduction}
\label{sec:intro}
Significant progress has been recently made in a variety of fields of information science using ideas centered around sparsity. Examples include compressed sensing~(Ref.~\cite{donoho2006compressed,Candes:2006robust}), matrix rank minimization (Ref.~\cite{recht2010guaranteed}), phase retrieval (Ref.~\cite{candes2013phase}) and robust  principal component analysis (Ref.~\cite{candes2011robust}), as well as many others. A key step in these examples is use of a variational formulation with a constraint or penaty term that is an $\ell_1$ or related norm. A limited set of extensions of sparsity techniques to physical sciences and partial differential equations (PDEs) have also appeared recently, including numerical solution of PDEs with multi-scale oscillatory solutions (Ref.~\cite{Schaeffer:2013sparse}) and efficient materials models derived from quantum mechanics calculations (Ref.~\cite{Ozolins2013compressive}). In all of these examples, sparsity is for the coefficients (i.e., only a small set of coefficients are nonzero) in a well-chosen set of modes (e.g., a basis or dictionary) for representation of the corresponding vectors or functions. In this paper, we propose a new use of sparsity techniques to produce ``compressed modes" - i.e., modes that are sparse and localized in space - for variational problems in mathematics and physics.

Our idea is motivated by the localized Wannier functions developed in solid state physics and quantum chemistry. We begin by reviewing the basic ideas for obtaining spatially localized solutions of the independent-particle Schr\"odinger's equation. For simplicity, we consider a finite system with $N$ electrons and neglect the electron spin. The ground state energy is given by $ E_0 =\sum_{j=1}^N \lambda_j$, where $\lambda_j$ are the eigenvalues of the Hamiltonian, $ \hat{H}=-\frac{1}{2} \Delta + V(\x)$, arranged in increasing order and satisfying $\hat{H} \phi_j = \lambda_j \phi_j$, with $\phi_j$ being the corresponding eigenfunctions. This can be recast as a variational problem requiring the minimization of the total energy subject to orthonormality conditions for wave functions:
\begin{equation}
E_0 = \min_{\Phi_N} \sum_{j=1}^N \langle \phi_j, \hat{H} \phi_j \rangle \quad \mbox{\text{s.t.}} \quad \langle \phi_j, \phi_k \rangle = \delta_{jk}.
\label{eqn:variationeigs_H}
\end{equation}
Here $\displaystyle \Phi_N = \{ \phi_j \}_{j=1}^N$ and $\displaystyle \langle \phi_j, \phi_k\rangle = \int_{\Omega} \phi^*_j(\x) \phi_k(\x)  \dx $ $(\Omega\subset\mathbb{R}^d)$. 

In most cases, the eigenfunctions $\phi_j$ are spatially extended and have infinite support, i.e., they are ``dense". This presents challenges for both computational efficiency (since the wave function orthogonalization requires $O(N^3)$ operations, dominating the computational effort for $N \approx 10^3$ electrons and above), and goes against physical intuition, which suggests that the screened correlations in condensed matter are short-ranged (Ref.~\cite{Prodan:1996}). 
It is well understood that the freedom to choose a particular unitary transformation (``subspace rotation") of the wave functions $\phi_j$ can be used to define a set of functions that span the eigenspace of $\hat{H}$, but are spatially localized or ``sparse".  Methods for obtaining such functions have been developed in solid state physics and quantum chemistry, where they are known as Wannier functions (Ref.~\cite{Wannier:1937}).

Mathematically, the Wannier functions are obtained as a linear combination of the eigenfunctions,
\begin{equation}
\label{eqn:WannierF}
W_j ({\bf x}) = \sum_k U_{jk} \phi_k ({\bf x}),
\end{equation}
where the subspace rotation matrix $U$ is unitary, $U^\dagger U=I$.  Currently, the most widely used approach to finding $W_j({\bf x})$ is the one proposed in Ref.~\cite{Marzari:1997} for calculating maximally localized Wannier functions (MLWFs). This approach starts with the pre-calculated eigenfunctions $\phi_j$ and determines $W_j$ by minimizing the second moment,
\begin{equation}
\label{eq:r2}
\langle \Delta {\x}^2_j \rangle = \langle W_{j}, ( \x - \langle {\x}_j \rangle)^2 W_{j} \rangle,
\end{equation}
where $\langle {\x}_j \rangle = \langle W_{j}, {\bf x} W_{j} \rangle$. More recently, a method weighted by higher degree polynomials is discussed in Ref.~\cite{E:2010PNAS}.  While this approach works reasonably well for simple systems, constructing optimally localized real-valued Wannier functions is often difficult because the minimization problem Eq.~(\ref{eq:r2}) is non-convex and requires a good starting point to converge to the global minimum. Another difficulty is that the resulting MLWFs need to be cut off ``by hand", which can result in significant numerical errors when the MLWFs are calculated ``on the fly" and their range is not known in advance. It would be highly desirable to devise an approach that does not require the calculation of the eigenfunctions and would converge to localized functions, while simultaneously providing a variational approximation to the total energy $E_0$. 

In this paper, we propose a novel method to create a set of  localized functions $\{\psi_i\}_{i=1}^N$, which we call compressed modes, 
such that $\sum_{j=1}^N \langle \psi_j, \hat{H} \psi_j \rangle$  approximates $E_0 = \sum_{j=1}^N \langle \phi_j, \hat{H} \phi_j \rangle$.
Our idea is inspired by the 
$\ell_1$ regularization used in compressive sensing. As a convex relaxation of $\ell_0$ regularization, $\ell_1$ regularization is commonly used for seeking sparse solutions for the underdetermined problem $A x = b$ (Ref.~\cite{Donoho:2003optimally,Candes:2006robust}). 
Motivated by advantages of the $\ell_1$ regularization for the sparsity in the discrete case, we propose a modification of the objective functional given by Eq.~(\ref{eqn:variationeigs_H}), which can immediately obtain functions with compact support and calculate approximate total energy ``in one shot", without the need to calculate eigenfunctions. This is accomplished by introducing an $L_1$ regularization of the wave functions:
\begin{equation}
\label{model:CMs}
E = \min_{\Psi_N} \sum_{j=1}^N  \left( \frac{1}{\mu}\left| \psi_j \right|_1 +  \langle \psi_j , \hat{H} \psi_j \rangle \right) \quad \mbox{\text{s.t.}} \quad \langle \psi_j, \psi_k \rangle = \delta_{jk},
\end{equation}
where $\Psi_N = \{\psi_j\}_{j=1}^N$ and the $L_1$ norm is defined as $\left| \psi_j \right|_1 = \int_\Omega | \psi_j | d\x$. For simplicity, we are requiring that the wave functions  $\psi_j$ are real; generalization to complex-valued wave functions, required to handle relativistic effects, is straightforward. The parameter $\mu$ controls the trade-off between sparsity and accuracy: larger values of $\mu$ will give solutions that better minimize the total energy at the expense of more extended wave functions, while a smaller $\mu$ will give highly localized wave functions at the expense of larger errors in the calculated ground state energy. Due to the properties of the $L_1$ term, the functions that minimize Eq.~(\ref{model:CMs}) will have compact support. In contrast to other approaches that use manually-imposed cutoff distances, the main advantage of our scheme is that one parameter $\mu$ controls both the physical accuracy and the spatial extent, while not requiring any physical intuition about the properties of the solution. In other words, the wave functions $\psi_j$ are nonzero only in those regions that are required to achieve a given accuracy for the total energy, and are zero everywhere else. Furthermore, due to the fact that exponentially localized Wannier functions are known to exist, the solution to Eq.~(\ref{model:CMs}) will provide a good approximation to the true total energy of the system (in fact, it converges to $E_0$ as $\mu^{-2}$).

There are two major contributions in this paper. First,  we propose localized compressed modes (CMs) using an $L_1$ regularized variational formula. In addition, we  propose a numerical algorithm to solve the proposed non-convex problem. Second, inspired by the idea of compressed modes, we further develop a new set of spatially localized orthonormal functions, compressed plane waves (CPWs), with multiresolution capabilities. Numerical algorithms - i.e., a fast compressed plane wave transform and a fast inverse compressed plane wave transform - are also developed using the new set of functions. 
%%%------------------------------------------------------------------------------------------------------------------

\section{Variational Model for Compressed Modes}
\label{sec:CMs}
\begin{figure}[h]
%\begin{minipage}{0.24\linewidth}
\centering
\includegraphics[width=.8\linewidth]{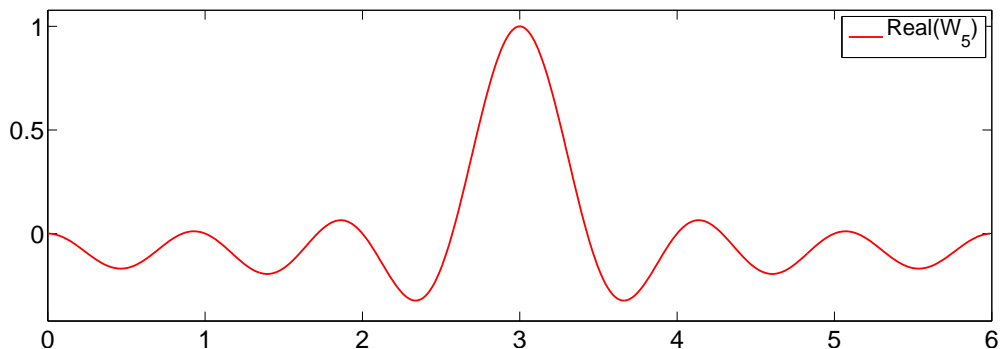}\\
%\end{minipage}\hfill
\label{fig:quasilocalfun}
\caption{A quasi-localized Wannier function for the 1D Laplace operator.}
\end{figure}

Consider a 1D free-electron case defined on $[0,L]$ with periodic boundary conditions. Namely, the Schr\"odinger operator is $\hat{H}_0 =-\frac{1}{2}\partial_{\x}^2$. It is clear that $\hat{H}_0$ has eigenfunctions $\frac{1}{\sqrt{L}}e^{i 2 \pi n x/L}$ with the corresponding eigenvalues $2(\pi n/L)^2, n = 0,\pm 1,\pm 2, \cdots$. With a unitary transformation $U_{mn} = \frac{1}{\sqrt{L}}e^{i 2\pi nm/L}$, one can construct quasi-localized orthonormal functions as
\begin{equation}
W_m = \frac{1}{L} \sum_{n} e^{i 2\pi n (x - m)/(L)}.
\end{equation}
Figure \ref{fig:quasilocalfun} illustrates the real part of one of the resulting quasi-localized functions obtained from the above unitary transformation. It is evident that the resulting $W_m$ are not even exponentially localized. 

\begin{figure}[h]
%\begin{minipage}{0.24\linewidth}
\centering
\includegraphics[width=.8\linewidth]{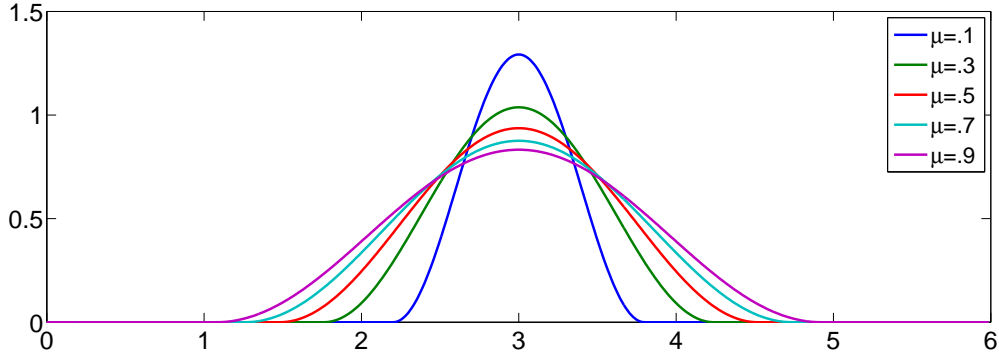}\\
%\end{minipage}\hfill
\label{fig:1Dpsi1_mu}
\caption{Theoretical $\psi_1$ in the 1D free electron model (Eq. (\ref{model:CM_psi_1}))  for different values of $\mu$.}
\end{figure}

As an example, we can analytically check that the $L_1$ regularization introduced in Eq. (\ref{model:CMs}) can localize the resulting functions. Let's again consider the 1D free-electron model defined on $[0,L]$ with the Schr\"odinger operator $\hat{H}_0 =-\frac{1}{2}\partial_{\x^2}$. Then the lowest mode satisfies 
\begin{eqnarray}
\psi_1 = \displaystyle\arg\min_{\psi} \frac{1}{\mu}\int_{\Omega} |\psi| \d x -\frac{1}{2} \int_{\Omega} \psi\partial_{\x^2}  \psi \d x 
%\hspace{2.2cm} \nonumber \\
            \quad \quad \text{s.t.} \quad \quad     \int_{\Omega} \psi(x)\psi(x)\d x =1.
\label{model:CM_psi_1} 
\end{eqnarray}
The solution of the above minimization problem will be a ``sparse'' solution, i.e., the Dirac delta function when $\mu \rightarrow 0$, and will approach the first 
eigenfunction of $\hat{H}$ when $\mu \rightarrow \infty$. Intuitively, we expect to be able to express the solution of Eq. (\ref{model:CM_psi_1}) as an approximation to a truncated diffusion of Dirac delta function via the Schr\"{o}dinger operator $-\frac{1}{2}\partial_{\x^2}$, which is a compactly supported function.
Indeed, the Euler-Lagrange equation corresponding to Eq. (\ref{model:CM_psi_1}) is
\begin{equation}
- \partial_{\x}^2 \psi_1 + \frac{1}{\mu} \text{sign}(\psi_1) = \lambda \psi_1 .
\label{pde:CM_psi_1}
 \end{equation}
 If we further assume that $\psi_1$ is symmetric around $x = L/2$,  the solution of [\ref{model:CM_psi_1}] is 
 \begin{equation}
\psi_1 = \left\{\begin{array}{cc}
                          \frac{1}{\lambda\mu}(1+ \cos (\sqrt{\lambda} (x-L/2))) & \text{if} \quad   |x-L/2| \leq l, \\
                          0& \text{if} \quad l \leq |x - L/2| \leq L,
                          \end{array}\right.
\end{equation}
where $  l = \pi/\sqrt{\lambda}$ and $\lambda =  (3\pi)^{2/5}\mu^{-4/5}$. Note that $\psi_1=\partial_{\x}\psi_1=0$ and $\partial_{\x}^2 \psi_1$ has a jump of $-\mu^{-1}$ at the boundary $x=l$ of the support of $\psi_1$, which are all consistent with Eq. (\ref{pde:CM_psi_1}). From this simple 1D example, it is clear that $L_1$ regularization can naturally truncate the resulting functions. Here $\psi_1$ has compact support $[L/2-l,L/2+l]$ if $\mu$ is small enough satisfying $l = \pi/\sqrt{\lambda}<L$. Moreover, we also observe that the smaller $\mu$ will provide a smaller region of compact support. Figure~\ref{fig:1Dpsi1_mu} shows $\psi_1$ for different values of $\mu$. 

This 1D example inspires us to consider the $L_1$ regularization of the wave functions proposed in Eq. (\ref{model:CMs}) for a general Schr\"odinger operator $\hat{H} = -\frac{1}{2}\Delta +V(\x)$ defined on $\Omega\subset \mathbb{R}^d$. 
\begin{definition}
We call $\Psi_N  = \{\psi_1,\cdots,\psi_N\}$ defined in the variational model \text{(\ref{model:CMs})} the first $N$ compressed modes (CMs) of the Schr\"{o}dinger operator $\hat{H}$. 
\end{definition}

By analogy with the localized Wannier functions described in the introduction, we expect that the CMs have compact support and can be expressed as orthonormal combinations of the eigenmodes of the original Schr\"{o}dinger operator. In other words, 
let $\Phi_M = \{\phi_1,\cdots,\phi_M\}$ be the first $M$ eigenfunctions of $\hat{H}$ satisfying
\begin{equation}
\hat{H} \Phi_M = \Phi_M \text{diag}(\lambda_1,\cdots,\lambda_M) \quad \& \quad  \int_{\Omega} \phi_j\phi_k \dx= \delta_{ij}. 
\end{equation}
We formulate the following conjecture to describe the completeness of the CMs.
\begin{conjecture}
Given $N\geq M$, consider the $N\times N$ matrix $\langle \Psi_N^T, \hat{H} \Psi_N\rangle$ with the $(j,k)-th$ entry defined by $\int_{\Omega} \psi_j H\psi_k \dx$ and let $(\sigma_1,\cdots,\sigma_M)$ be its first $M$ eigenvalues; then
\begin{eqnarray}
\displaystyle\lim_{\mu \rightarrow\infty}\sum_{j = 1}^M (\sigma_j - \lambda_j)^2 = 0  \quad \text{and} \quad 
\lim_{N \rightarrow\infty}\sum_{j = 1}^M (\sigma_j - \lambda_j)^2 = 0 .
\end{eqnarray} \label{thm:completeness}
\end{conjecture}

%%%-------------------------------------------------------------------------------------------------------------------
\subsection{Numerical algorithms}

\begin{figure}[h]
\centering
\begin{minipage}{0.49\linewidth}
\centering  
\includegraphics[width=1\linewidth]{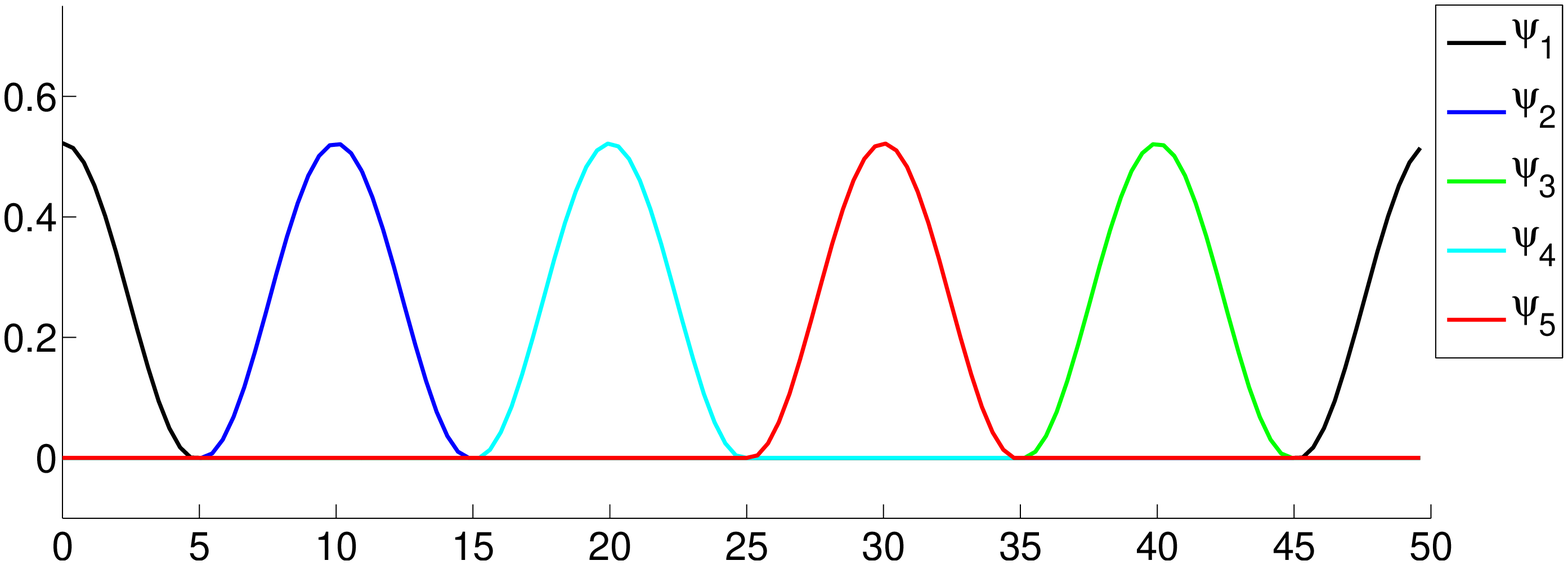}\\
\centering  $\mu = 30 $ 
\end{minipage}\hfill
\begin{minipage}{0.49\linewidth}
\centering
\includegraphics[width=1\linewidth]{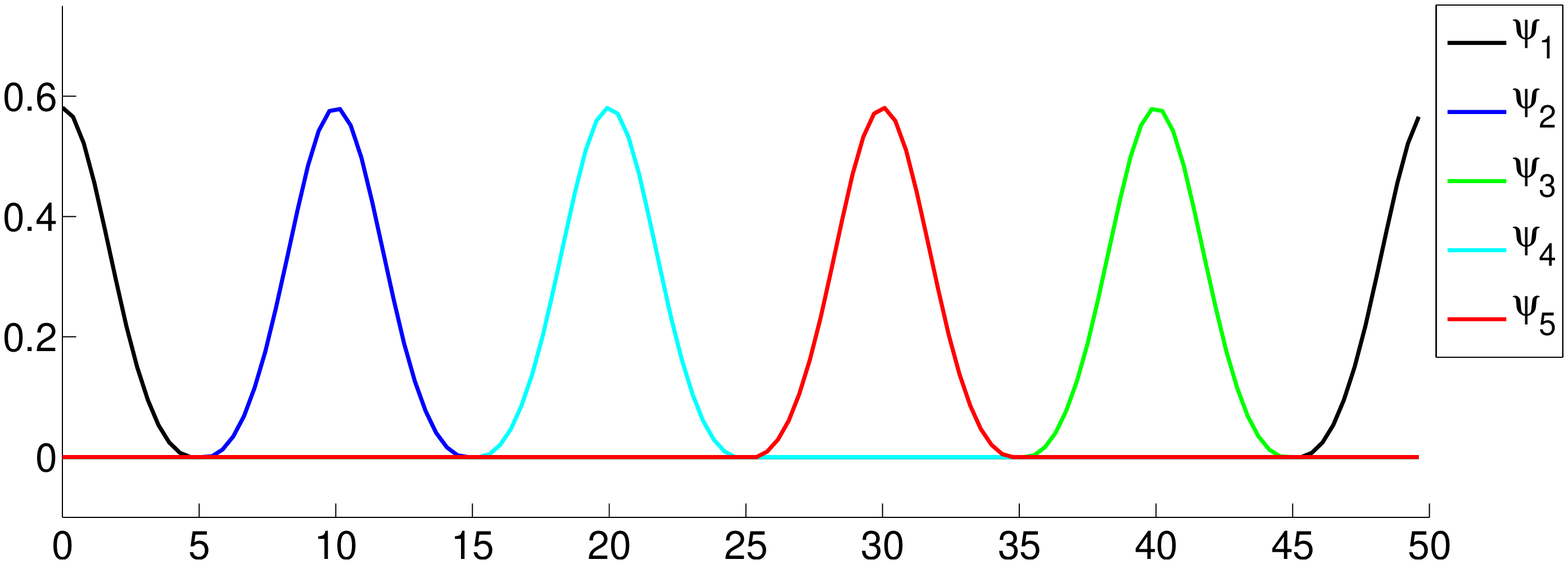}\\
\centering $\mu = 50$
\end{minipage}\hfill\\
\begin{minipage}{0.49\linewidth}
\centering
\includegraphics[width=1\linewidth]{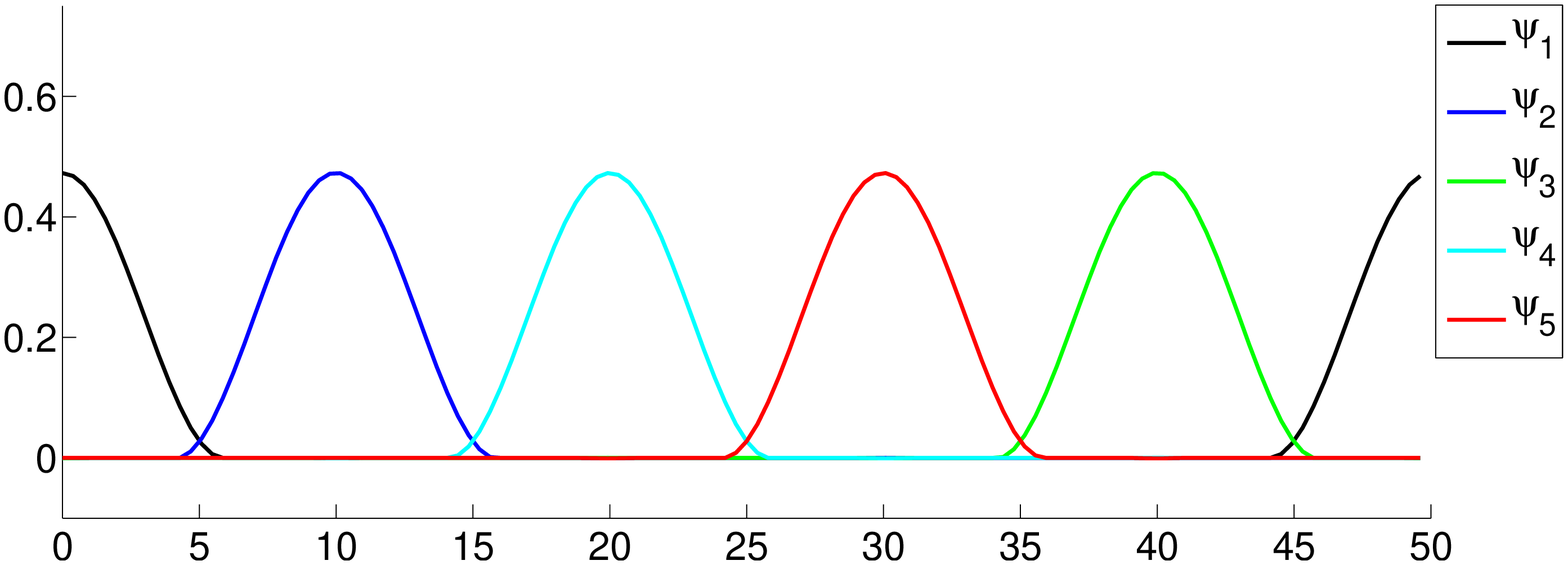}\\
\centering $\mu = 50$
\end{minipage}\hfill
\begin{minipage}{0.49\linewidth}
\centering
\includegraphics[width=1\linewidth]{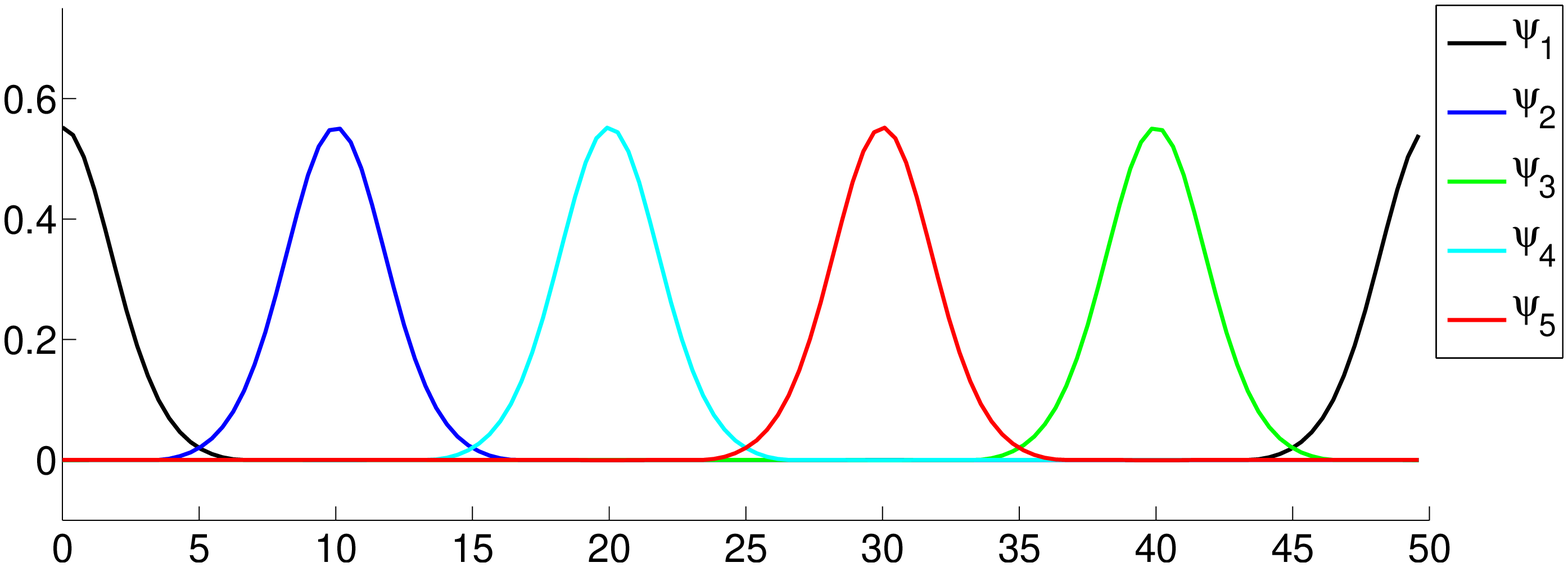}\\
\centering $\mu = 300$
\end{minipage}\hfill\\
\begin{minipage}{0.49\linewidth}
\centering
\includegraphics[width=1\linewidth]{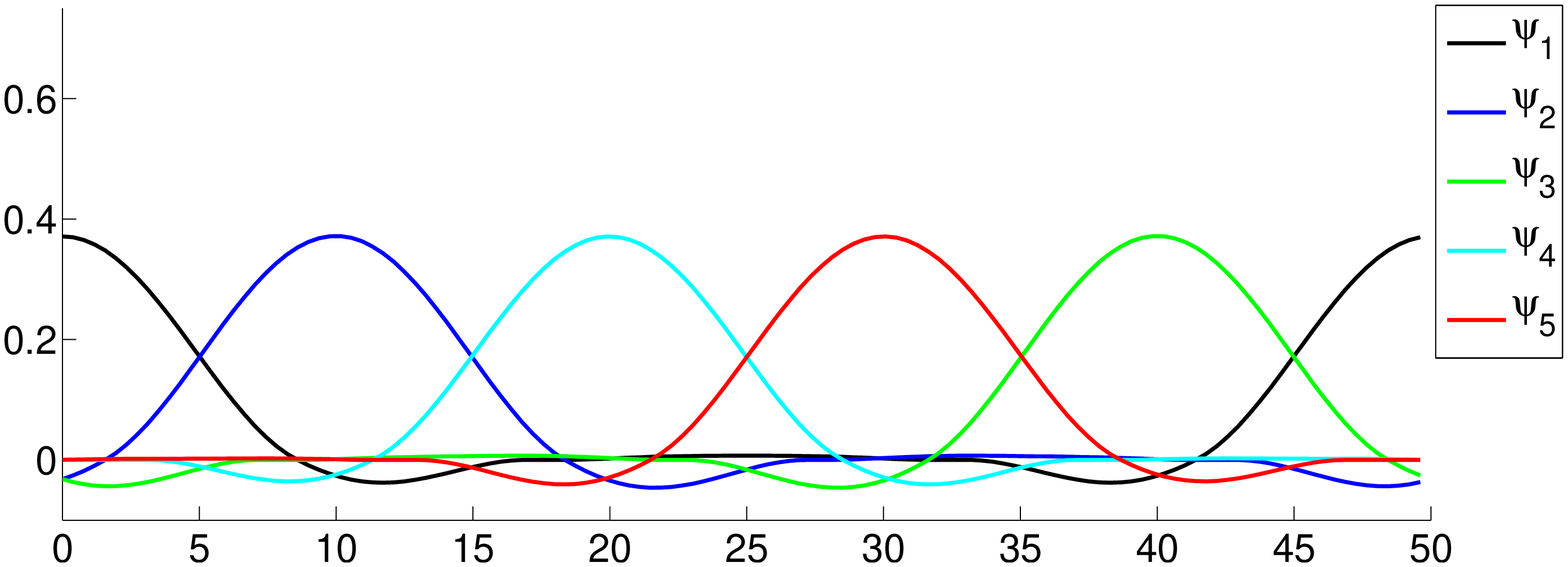}\\
\centering $\mu = 500$
\end{minipage}\hfill
\begin{minipage}{0.49\linewidth}
\centering 
\includegraphics[width=1\linewidth]{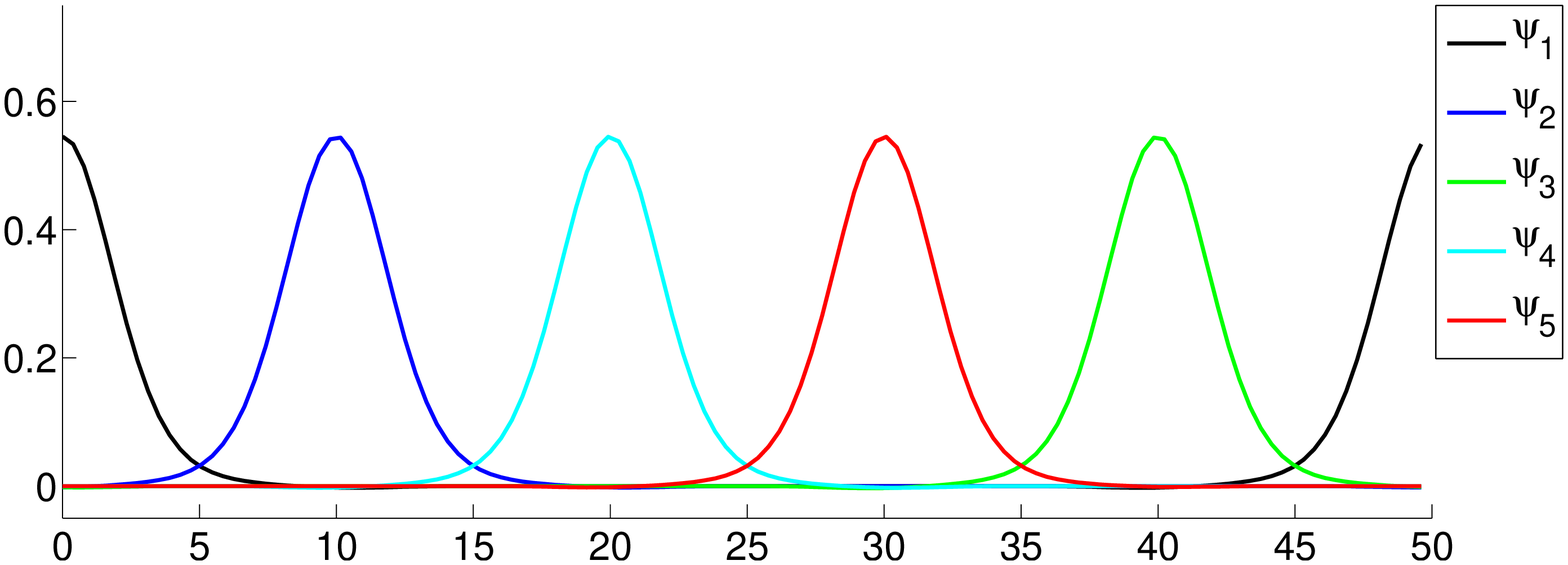}\\
\centering $\mu = 5000$
\end{minipage}\hfill\\
\caption{Computation results of CMs with different values of $\mu$. The first column: the first 5 CMs of the 1D free-electron model. The second column: the first 5 CMs of the 1D Kronig-Penny model. 
}
\label{fig:CMs}
\end{figure}

To numerically compute the proposed CMs, we consider the system  on a domain $D = [0,L]^d \subset \mathbb{R}^d$ with periodic boundary conditions, and  discretize the domain $D$ with $n$ equally spaced nodes in each direction. Then 
the variational formula Eq. (\ref{model:CMs}) for the first $N$-CMs can be reformulated and discretized as follows: 
 \begin{eqnarray}
\Psi_N = \min_{\Psi\in\mathbb{R}^{n\times N}} \frac{1}{\mu}|\Psi|  + \Tr \langle\Psi^T \hat{H} \Psi \rangle  \quad \text{s.t.} \quad \Psi^T\Psi = I 
\label{eqn:CEM_discretization},
\end{eqnarray}
in which $|\Psi|$ is the $\ell_1$ norm of the matrix $\Psi$.

We  solve this optimization problem using the algorithm of splitting orthogonality constraint (SOC) proposed in Ref.~\cite{Lai:2013JSC}.
By introducing auxiliary variables $Q = \Psi$ and $P= \Psi$,  the above problem Eq. (\ref{model:CM_psi_1}) is equivalent to the following constrained problem:
 \begin{eqnarray}
\min_{\Psi,P,Q}\frac{1}{\mu} |Q| + \Tr \langle\Psi^T \hat{H} \Psi \rangle 
\quad \text{s.t.}  \quad  Q=\Psi,  P = \Psi,  P^TP = I,
\label{eqn:CEM_splitting}
\end{eqnarray}
which can be solved by the SOC algorithm based on split Bregman iteration (Refs.~\cite{Osher:2005,Yin:2008bregman,Goldstein:2009split}).

\begin{algorithm}
\label{alg:CM_SOC}
Initialize $\Psi_N^0 = P^0 = Q^0 , b^0 = B^0 = 0$.

\While{``not converged"}{
\begin{enumerate}
\item $\displaystyle \Psi_N^k = \arg\min_{\Psi} \Tr \langle\Psi^T \hat{H} \Psi \rangle+ \frac{\lambda}{2} \|\Psi - Q^{k-1}
     + b^{k-1}\|^2_F+ \frac{r}{2} \|\Psi - P^{k-1} + B^{k-1}\|_F^2 $.
\item $\displaystyle Q^k = \arg\min_{Q} \frac{1}{\mu} |Q| + \frac{\lambda}{2} \| \Psi_N^k - Q + b^{k-1}\|_F^2  $.
\item $P^k = \displaystyle\arg\min_{P} \frac{r}{2} \|\Psi_N^k - P +B^{k-1}\|_F^2 \quad  \text{\rm s.t.} \quad P^TP =I $.
\item $b^{k} = b^{k-1} +  \Psi_N^k - Q^k$.
\item $B^{k} = B^{k-1} +  \Psi_N^k - P^k$.
\end{enumerate}
}\end{algorithm}
The above submimization problems can be easily solved as follows:
\begin{eqnarray}
\displaystyle \left( 2\hat{H} + \lambda + r\right ) \Psi_N^k &=& r(P^{k-1} - B^{k-1}) + \lambda(Q^{k-1} - b^{k-1}). \label{eqn:schrodinger}\\
\displaystyle Q^k  &=& \mbox{\text{Shrink}}(\Psi_N^k + b^{k-1},1/(\lambda\mu))  \label{eqn:shrinkage}\\
 \displaystyle P^k &=& (\Psi_N^k + B^{k-1})U\Lambda^{-1/2}S^T  \label{eqn:quadorthogonality}
\end{eqnarray}
where $U\Lambda S^T= \mbox{svd}\left((\Psi^k + B^{k-1})^T(\Psi^k + B^{k-1})\right)$  and the ``Shrink" (or soft-threshholding) operator is defined as $\text{Shrink}(u,\delta)=\text{sgn}(u) \max(0,|u|-\delta)$. Since the matrix $ 2\hat{H} + \lambda + r$ in Eq. (\ref{eqn:schrodinger}) is sparse and can be made to be positive definite,  a few iterations of either Gauss-Seidel or conjugate gradient will be good enough in practice. Thus, Eq. (\ref{eqn:schrodinger}) and (\ref{eqn:shrinkage}) can be solved very efficiently with long operation counts linearly dependent on $N$. The only time consuming part in our algorithm is Eq. 
[\ref{eqn:quadorthogonality}], which involves an SVD factorization and can be straightforwardly solved with an $O(N^3)$ algorithm. For a moderate size of $N$, the proposed algorithm can solve the problem efficiently. 

For a large number of modes $N$, one possible approach to speed up the computation is to use GPU based parallel computation to perform the SVD factorization. 
Here, we propose another method to speed up the proposed algorithm due to the special structure of the solution. 
We find that each of the resulting functions  $\psi_1,\cdots,\psi_N$ has compact support, so that the support of each $\psi_i$ overlaps with only a finite number of its
neighbors. This allows us to replace to the full orthogonality constraint $ \Psi^T\Psi = I$ by a system of banded orthogonality constraints.
\begin{equation}
\displaystyle \int \psi_j\psi_k \dx= \delta_{jk}, 
\quad \left\{\begin{array}{c}
j = 1,\cdots, N. \\
k = j, j\pm 1, \ldots j\pm p
\end{array}\right.
\end{equation}
where $p$ is the band width. Thus, the $O(N^3)$ algorithm for SVD factorization in Eq. (\ref{eqn:quadorthogonality}) can be replaced by $N$ factorizations of $2p\times 2p$ matrices, which is an algorithm with $8p^3O(N)$ long operations.

Note that the discretized variational principle in Eq.~(\ref{eqn:CEM_discretization}) is related to sparse principal component analysis (SPCA) (Ref.~\cite{AspremontSPCA2007,Qi:2013sparse}).  SPCA, however,  does not involve an underlying continuum variational principle and the sparse principal components are not localized, since the component number does not correspond to a continuum variable.

%%%-------------------------------------------------------------------------------------------------------------------
\subsection{Numerical results} 
In quantum mechanics, the free electron model is a simple model for the behavior of valence electrons in a metallic solid with weak atomic pseudopotentials, where the potential function in the Schr\"{o}dinger operator is simply set as zero. The Kronig--Penney (KP) model (Ref.~\cite{Kronig:1931quantum}) is a simple, idealized quantum-mechanical system that describes the movement of independent electrons in a one-dimensional crystal, where the potential function $V({\bf x})$ consists of a periodic array of potential wells described by $V(x) = -V_0\sum_{j=1}^{Nel}e^{-\frac{(x - x_j)^2}{2\delta^2}}$. We choose $Nel = 5, V_0 = 1, \delta = 3$ and $x_j = 10j$ in our 1D KP model discussed below.Ê

\begin{figure}[h]
\centering
\begin{minipage}{0.49\linewidth}
\centering
\includegraphics[width=1\linewidth]{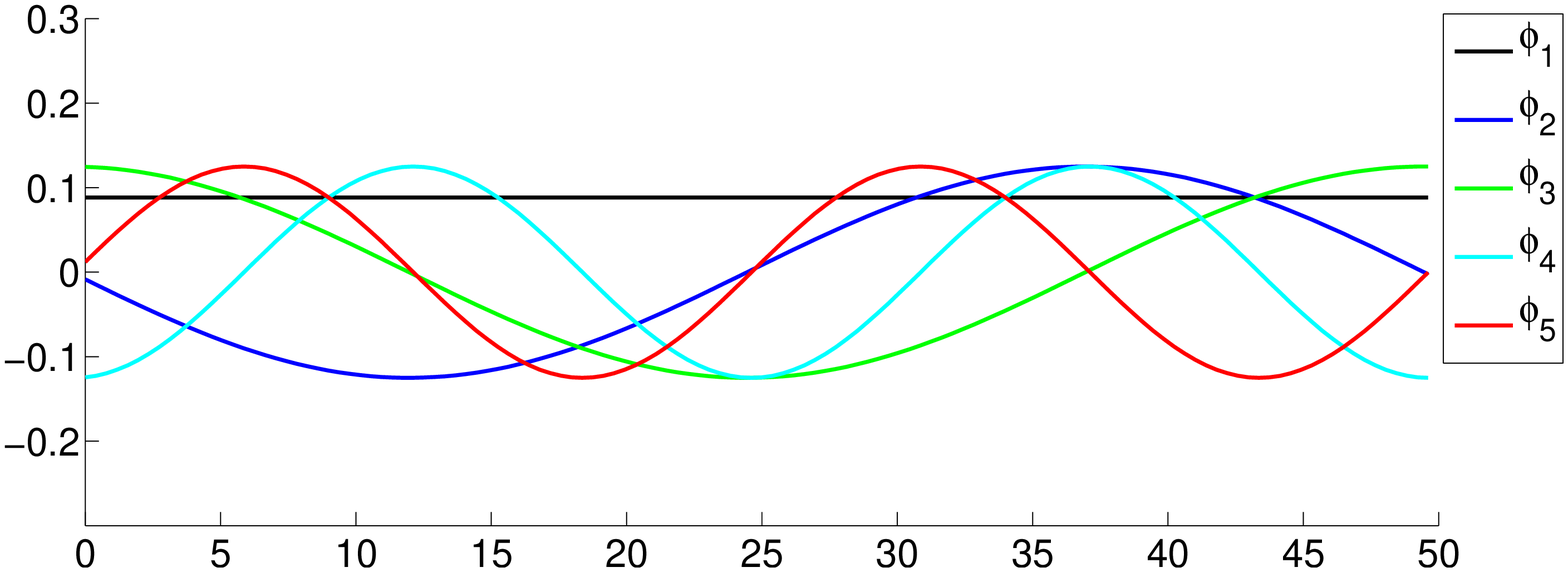}\\
\end{minipage}\hfill
\begin{minipage}{0.49\linewidth}
\centering
\includegraphics[width=1\linewidth]{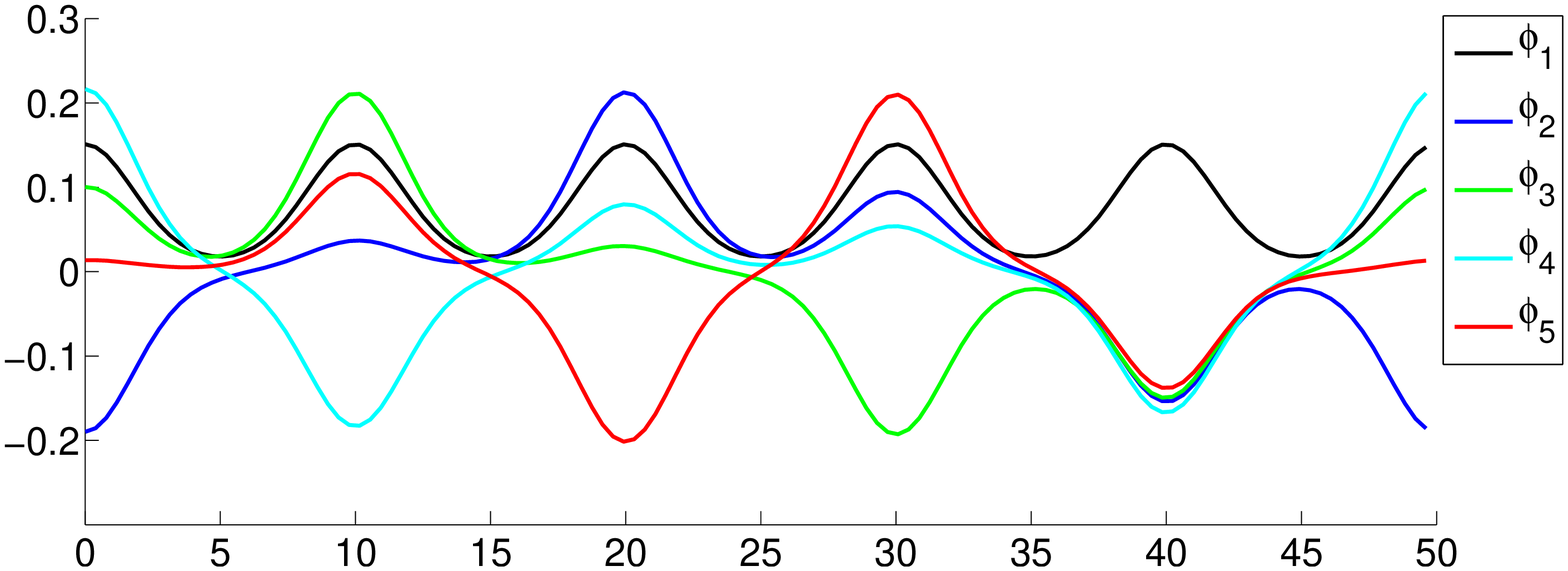}\\
\end{minipage}\hfill
\caption{The first 5 egienfunctions of the Schr\"{o}dinger operator $\hat{H}$ in the free-electron model (left) and the KP model (right), respectively.}
\label{fig:FE_KP_eigs}
\end{figure}

In our experiments, we choose $\Omega = [0, 50]$ 
and  discretize $\Omega$ with 128 equally spaced nodes.  The proposed variational model Eq. (\ref{model:CMs}) is solved using  algorithm~\ref{alg:CM_SOC}, where parameters are chosen as $\lambda = \mu N/20$ and $r = \mu N/5$. We report the computational results of the first 5 CMs of the 1D free-electron model (the first column) and the 1D KP model (the second column) 
in Figure \ref{fig:CMs}, where we use 5 different colors to differentiate these CMs. To compare all results more clearly, we use the same initial input for different values of $\mu$ in the free-electron model and the 1D KP model. We flip the CMs if necessary such that most values of CMs on their support are positive, since sign ambiguities do not affect minimal values of the objective function in Eq. (\ref{model:CMs}).
 For comparison, Figure~\ref{fig:FE_KP_eigs} plots the first 5 eigenfunctions of the Schr\"{o}dinger operator used in the free-electron model and KP model. It is clear that all these eigenfunctions are spatially extended without any compact support. However, as we can observe from Figure \ref{fig:CMs}, the proposed variational model does provide a series of compactly supported functions.
Furthermore, all numerical results in Figure~\ref{fig:CMs} clearly show the dependence of compact support size via $\mu$ as it can be prescribed by the variational formula Eq. (\ref{model:CMs}). In other words,  the model with smaller $\mu$ will create CMs with smaller compact support, and the model with larger $\mu$ will create CMs with larger compact support. In addition, we can also find that the resulting compressed modes are not interacting for small $\mu$ (the first row of Figure~\ref{fig:CMs}). By increasing $\mu$ to a moderate value, the modes start to interact each other with a small amount of overlap (the second row of Figure~\ref{fig:CMs}). Significant overlap can be observed using a big value of $\mu$ (the second row of Figure~\ref{fig:CMs}).

\begin{figure}[h]
\begin{minipage}{0.33\linewidth}
\centering $\mu = 10, M = 50, N = 50$
\end{minipage}\hfill
\begin{minipage}{0.33\linewidth}
\centering $\mu = 10, M = 50, N = 60$
\end{minipage}\hfill
%\begin{minipage}{0.24\linewidth}
%\centering $\mu = 10, M = 50, N = 100$
%\end{minipage}\hfill
\begin{minipage}{0.33\linewidth}
\centering $\mu = 10, M = 50, N = 128 $
\end{minipage}\hfill\\
\begin{minipage}{0.33\linewidth}
\centering
\includegraphics[width=1\linewidth]{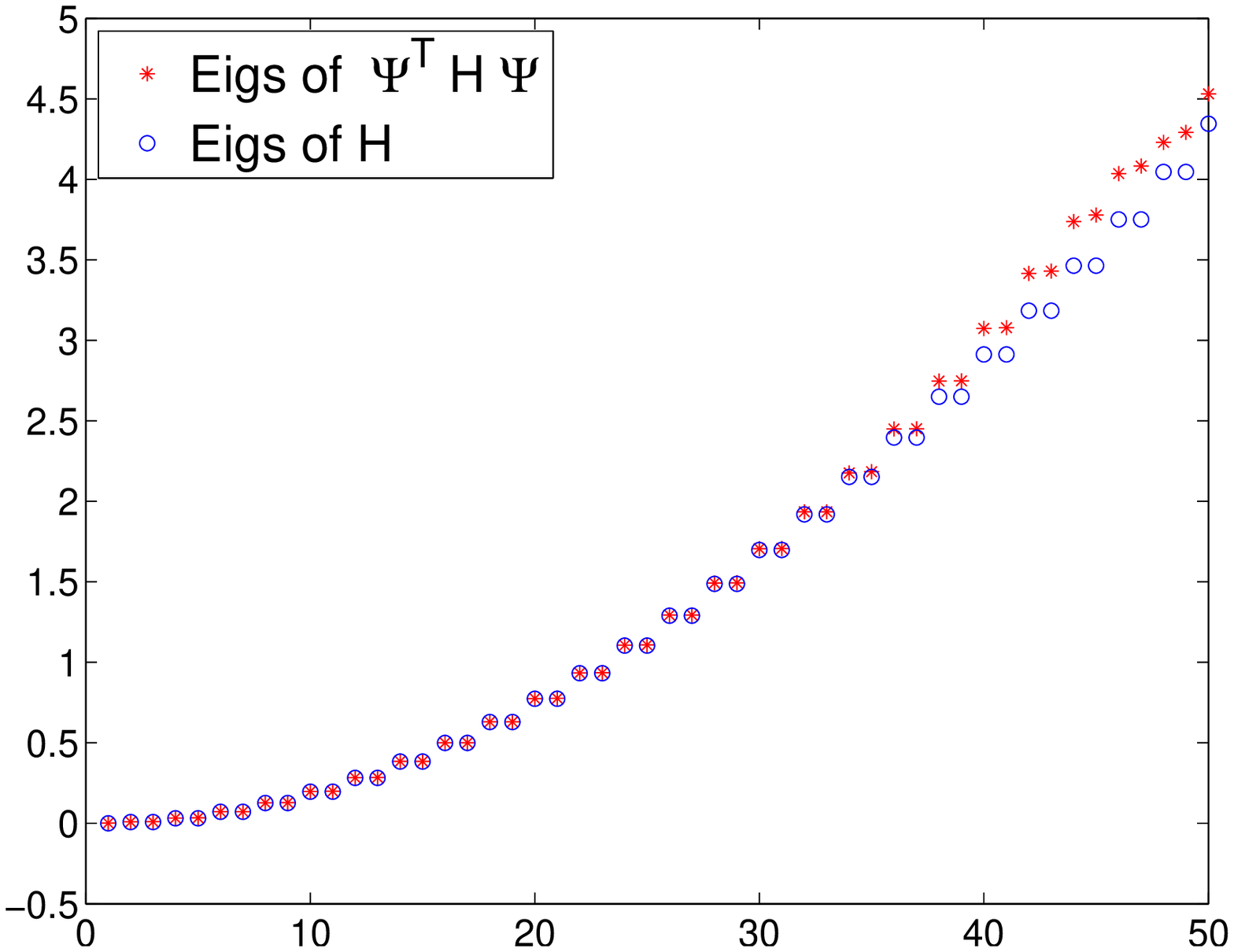}
\end{minipage}\hfill
\begin{minipage}{0.33\linewidth}
\centering
\includegraphics[width=1\linewidth]{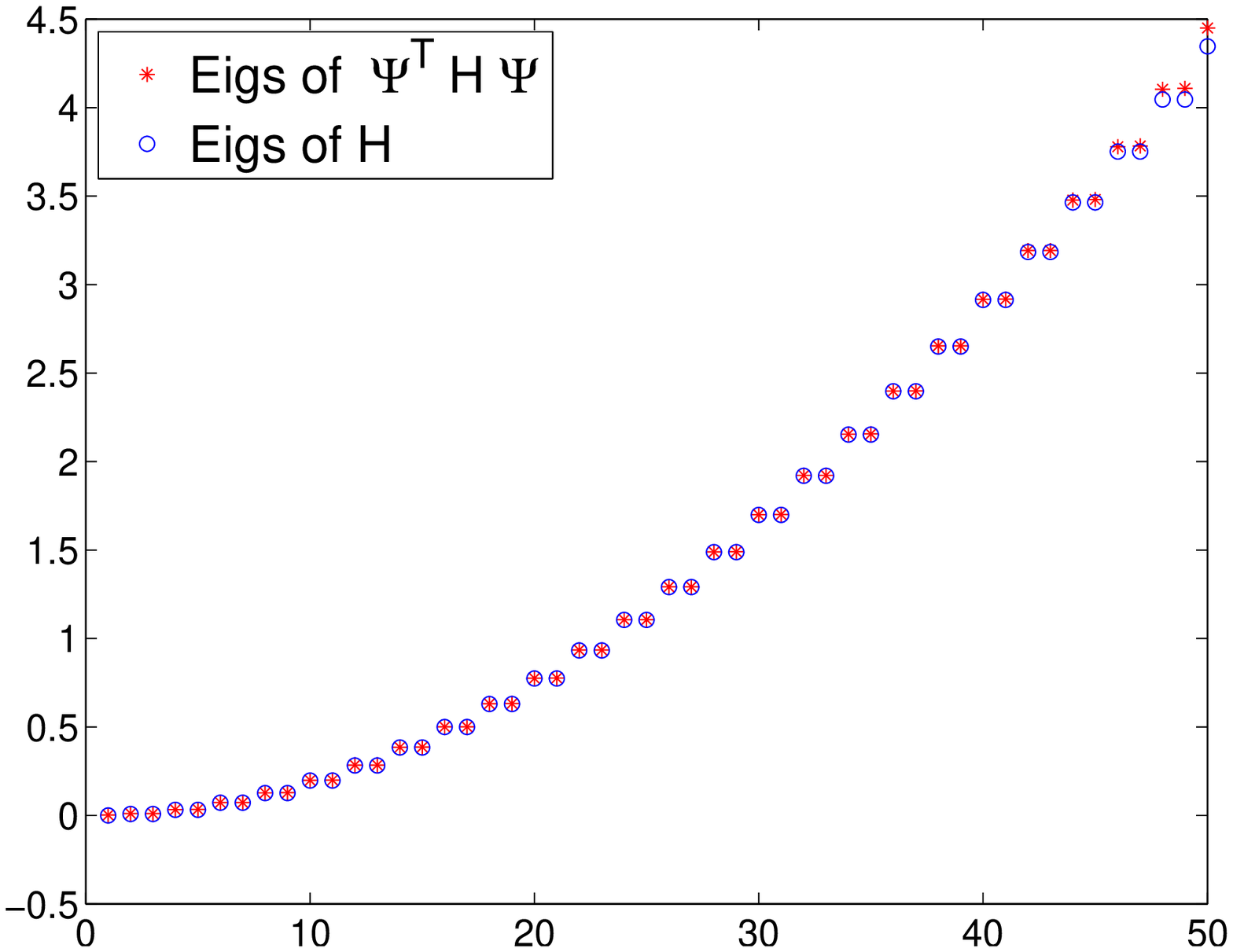}
\end{minipage}\hfill
%\begin{minipage}{0.24\linewidth}
%\centering
%\includegraphics[width=1\linewidth]{CMs_1D_FE_ErrorEigs_mu10_N100.eps}
%\end{minipage}\hfill
\begin{minipage}{0.33\linewidth}
\centering
\includegraphics[width=1\linewidth]{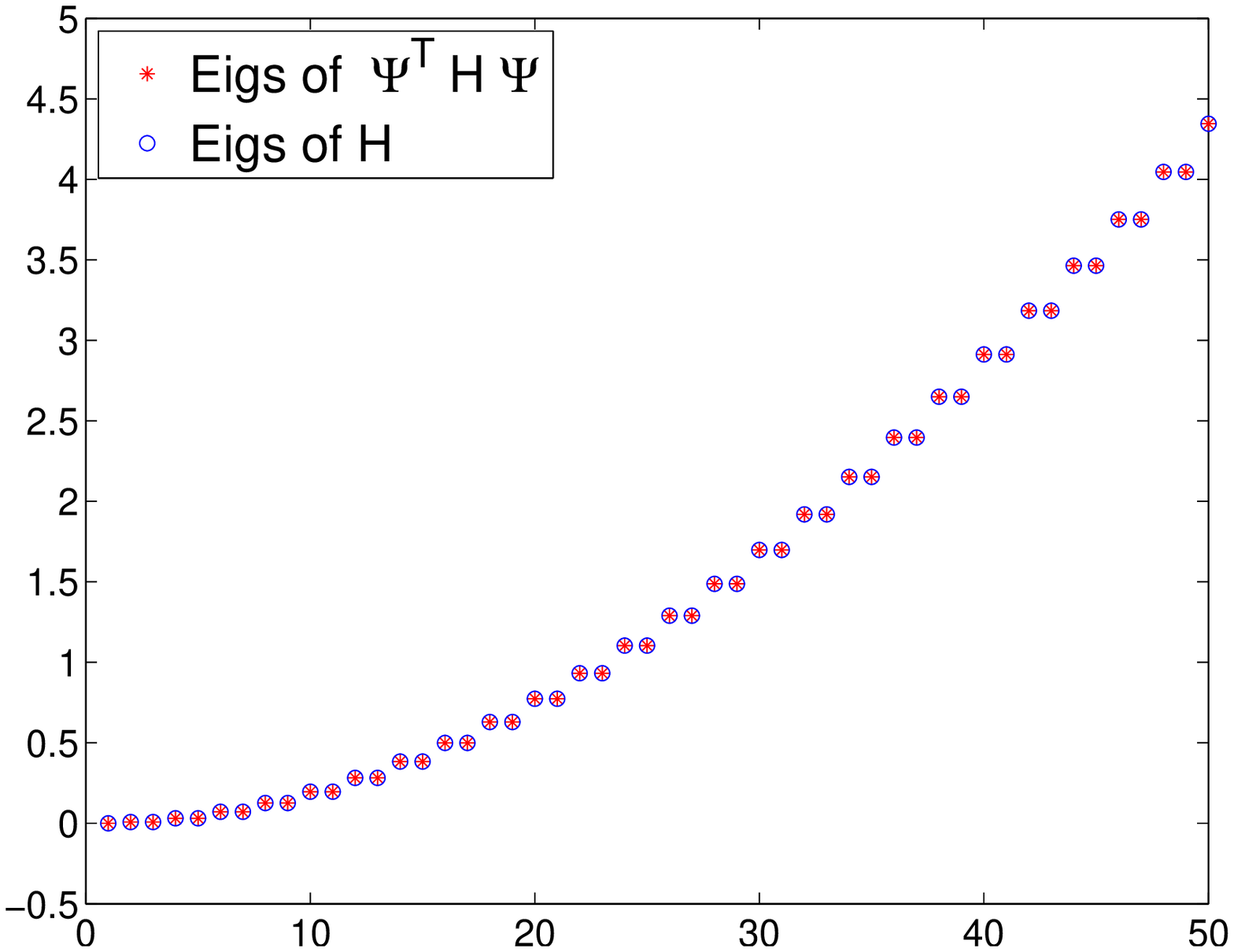}
\end{minipage}\hfill\\
\begin{minipage}{0.33\linewidth}
\centering
\includegraphics[width=1\linewidth]{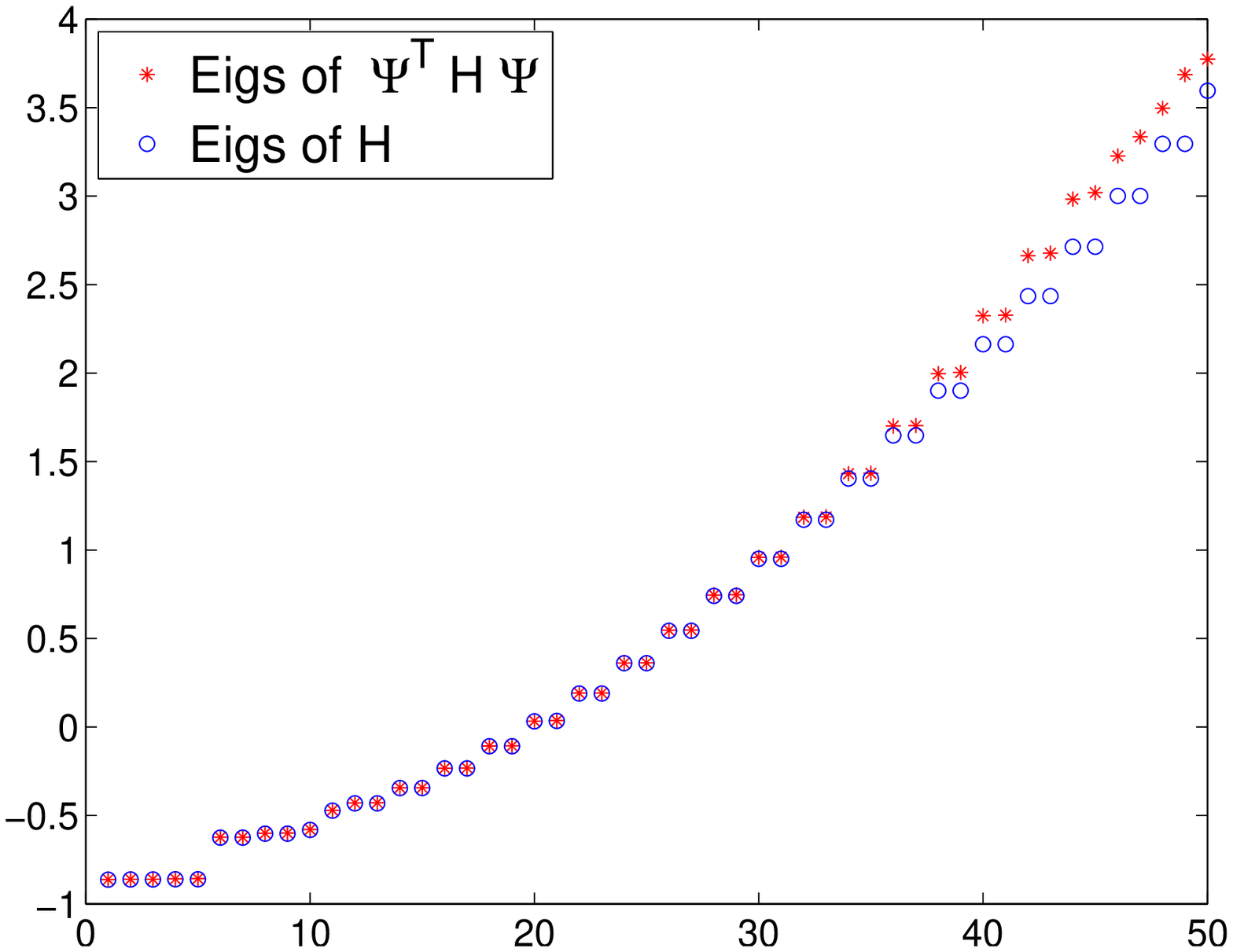}
\end{minipage}\hfill
\begin{minipage}{0.33\linewidth}
\centering
\includegraphics[width=1\linewidth]{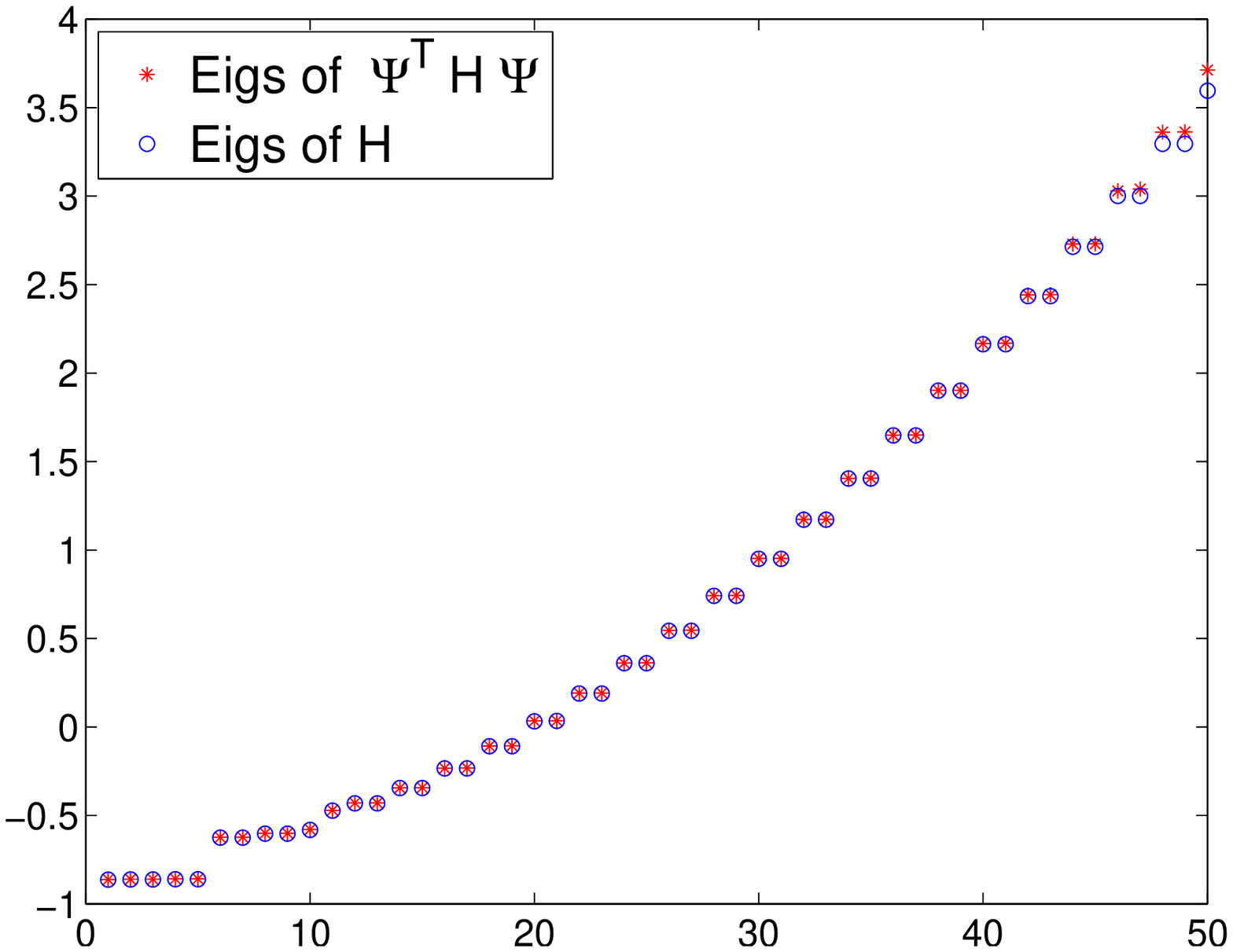}
\end{minipage}\hfill
%\begin{minipage}{0.24\linewidth}
%\centering
%\includegraphics[width=1\linewidth]{CMs_1D_KP_ErrorEigs_mu10_N100.eps}
%\end{minipage}\hfill
\begin{minipage}{0.33\linewidth}
\centering
\includegraphics[width=1\linewidth]{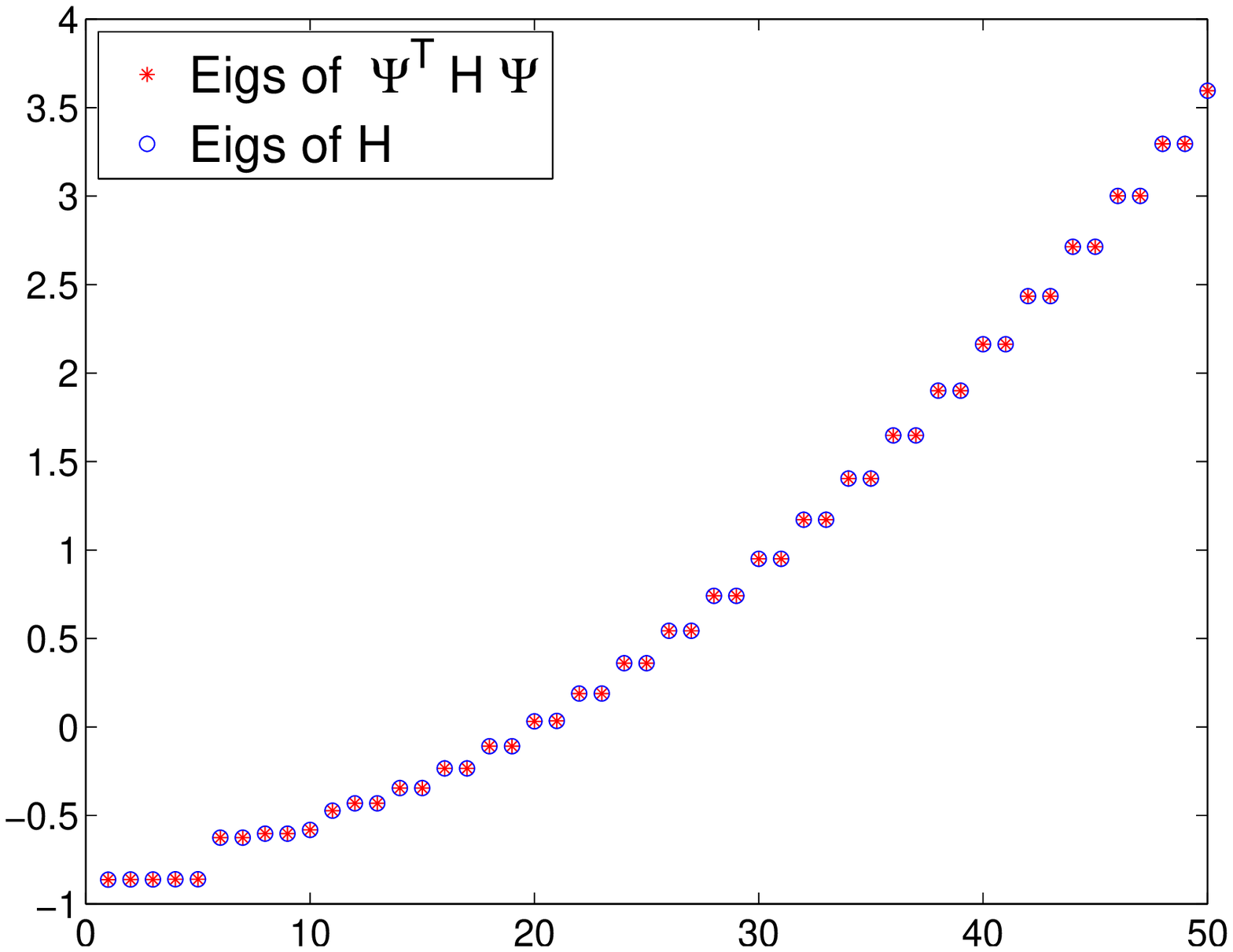}
\end{minipage}\hfill\\
\caption{Comparisons of the first 50 eigenvalues of the 1D free electron model (the first row) and the 1D Kronig-Penney model (the second row). }
\label{fig:EigsErrorComparisons1}
\end{figure}

\begin{figure}[htp]
\centering
\begin{minipage}{0.49\linewidth}
\centering
\includegraphics[width=.92\linewidth]{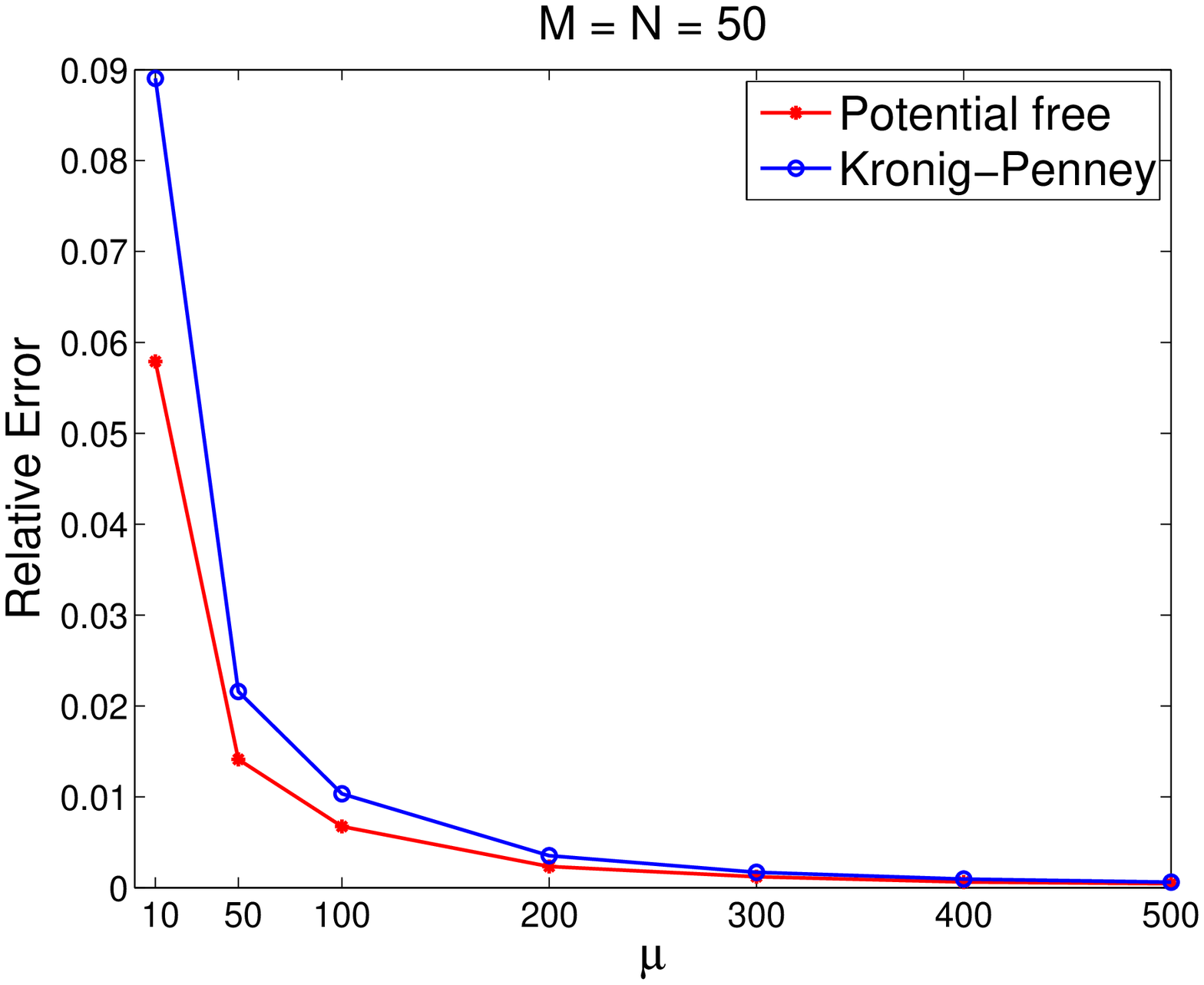}
\end{minipage}\hfill
\begin{minipage}{0.49\linewidth}
\centering
\includegraphics[width=.9\linewidth]{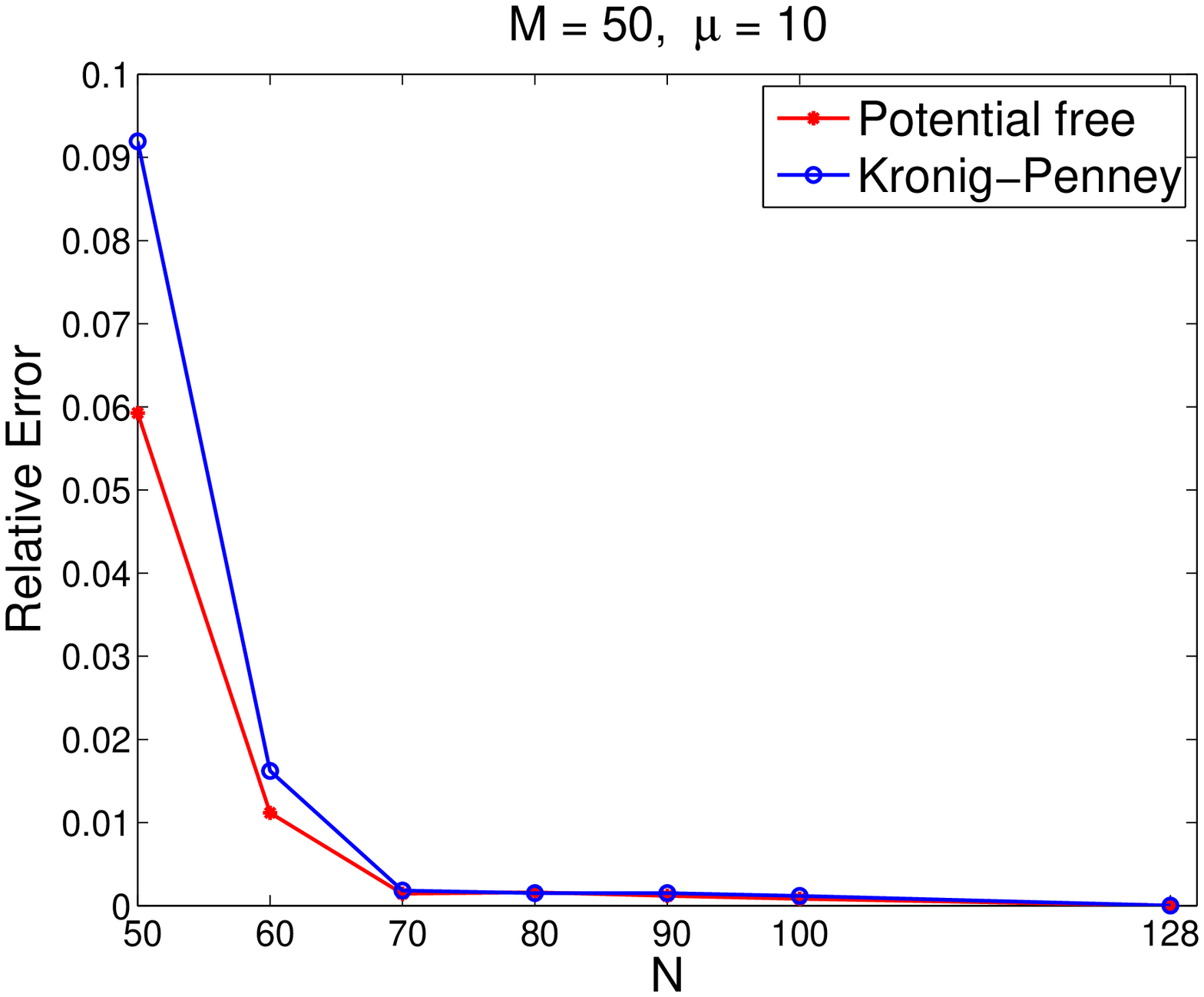}
\end{minipage}\hfill\\
\caption{Relative eigenvalue error of the 1D free-electron model (red dots) and 1D KP model (blue circles). Left: relation of the relative error via different values of $\mu$ for fixed $M = N = 50$. Right: relation of the relative error via different values of $N$ for fixed $\mu = 10$ and $M = 50$.}
\label{fig:EigsErrorComparisons2}
\end{figure}

Inspired by Wannier functions, which are compactly supported and can be represented as orthonormal combinations of the eigenmodes of the Schr\"{o}dinger operator, we also compare the first $M$ eigenvalues $(\sigma_1,\cdots,\sigma_M)$ of the matrix $\langle \Psi_N^TH \Psi_N\rangle$ obtained by the 1D KP model and 1D free-electron model with the first $M$ eigenvalues $(\lambda_1,\cdots,\lambda_M)$ of the corresponding Schr\"{o}dinger operators. Figure~\ref{fig:EigsErrorComparisons1} illustrates the comparisons with a relative small value $\mu=10$. We can clearly see that $\{\sigma_i\}$ gradually converges to $\{\lambda_i\}$ as $N$ goes to the number of the nodes in our discretization. In addition, we also plot the relative error $E = \sum_{i=1}^M(\sigma_i - \lambda_i)^2/\sum_{i=1}^M(\lambda_i)^2$ in Figure~\ref{fig:EigsErrorComparisons2}. As we speculated in conjecture \ref{thm:completeness}, the relative error will converge to zero as $\mu\rightarrow\infty$ for fixed $M = N = 50$, which is illustrated in the left image of Figure~\ref{fig:EigsErrorComparisons2}. The relative error will also converge to zero as $N\rightarrow\infty$ for fixed $\mu = 10$ and $M = 50$, which is illustrated in the right image of Figure~\ref{fig:EigsErrorComparisons2}. 

All the above numerical experiments validate that the proposed CMs do provide a series of compactly supported orthonormal functions. Moreover, the CMs are found to approximately span the eigenspace  of the Schr\"{o}dinger operator (i.e., the space of linear combinations of its first few lowest eigenmodes). 

%%%-------------------------------------------------------------------------------------------------------------------
\section{Compressed Plane Waves (CPWs)}
\label{sec:CPW}
Plane waves/Fourier basis functions have countless important applications in science and engineering. Plane waves are all global functions, while localized basis functions are remarkably useful in fields such as signal processing, image science and data science. The lack of localized properties of plane waves has been complemented by wavelet functions (Ref.~\cite{Daubechies:1992ten}), which have been tremendously successful in many applications. However, canonical wavelet functions usually can only be defined on regular domains in $\mathbb{R}^d$ by tensor products of wavelet functions in 1D. More recently, a diffusion wavelet is proposed in Ref. \cite{Coifman2006diffusion}, which can be considered on general domains, manifolds and graphs. It essentially can be viewed as diffusion of delta functions where multiresolution can be obtained by choosing different diffusion times $t$.  This method is completely based on a diffusion process without considering any sparsity-inducing variational approach. Inspired by the $L_1$ regularized variational method for generating compressed modes proposed in the previous section, we would like propose a new method of generating localized orthonormal functions adapted to a given differential operator. 

We consider a simple elliptic operator $-\Delta/2$ defined in $\mathbb{R}^d$, which also corresponds to the free electron case with the Schr\"{o}dinger operator $\hat{H}_0 = -\Delta/2$. Let $\w = (w_1,\cdots,w_d)\in\mathbb{R}_{+}^d$ be a basis of the $d$-dimensional lattice $\Gamma_{\w} = \{\j\w:=(j_1w_1,\cdots,j_dw_d)\hspace{0.1cm} |\hspace{0.1cm}\j=(j_1,\cdots,j_d)\in\mathbb{Z}^d\}$.  Similar to the variational model for the CMs, we introduce a set of localized orthonormal  modes $\{\psi^n\}_{n=1}^\infty$ recursively defined as follows:
\begin{eqnarray}
%\left\{\begin{array}{c}
\psi^1 = \displaystyle\arg\min_{\psi}\frac{1}{\mu} \int_{\Omega} |\psi| \dx + \int_{\Omega}  \psi \hat{H}_0\psi \dx %\hspace{2.2cm} \nonumber \\
            \quad \text{s.t.} \quad
            \int_{\Omega}  \psi(\x)\psi(\x - \j\w)\dx = \delta_{\j0}, \quad  \j\in\mathbb{Z}^d.  \hspace{1.2cm}  \label{eqn:CPW1}
\end{eqnarray}
The higher modes can be recursively defined as:
\begin{eqnarray}
\psi^n =\displaystyle \arg\min_{\psi} \frac{1}{\mu} \int_{\Omega}  |\psi| \dx + \int_{\Omega}  \psi \hat{H}_0\psi \dx %\hspace{2.2cm} \nonumber \\
        \quad  \text{s.t.}  \quad
           \left\{\begin{array}{cc}\displaystyle \int_{\Omega}  \psi(\x)\psi(\x - \j\w)\dx = \delta_{\j0}, &  \j\in\mathbb{Z}^d  \vspace{0.2cm}\\ 
           \displaystyle\int_{\Omega} \psi(\x)\psi^i(\x - \j\w)\dx = 0,  & i = 1,\cdots, n-1
           \end{array}\right.  \label{eqn:CPWn}
%\end{array}\right. \\
\end{eqnarray}
Here, the parameters $\mu$ and $\w$ are given. Recall that our goal is to construct a set of localized orthonormal functions which can
provide a new tool to analyze functions. With the help of the localized orthonormal  modes $\{\psi^n\}_{n=1}^\infty$, we can construct a new set of orthonormal functions described as follows.

\begin{definition}
We define $\{\b^n_\j(\x) = \psi^n(\x - \j\w) \}^{n = 1,2,3\cdots}_{\j\in\mathbb{Z}^d}$.
We  call $\{\psi^n\}_n$ the basic compressed plane waves (BCPWs) and call $\{\b^n_\j\}_{n,\j}$ the compressed plane waves (CPWs). 
\end{definition}

Ideally, we hope that the CPWs are complete and form an orthonormal basis. Inspired by our numerical experiments (see Fig. \ref{fig:CPW_spectraldensity}), we formulate  the following conjecture for the completeness of CPWs. This will be studied in our future work.

\begin{conjecture}[Completeness of CPWs]
There exists a constant $\mu_0$ such that the set of orthonormal functions $\{\b^n_\j\}_{n,\j}$ generated from Eqs. (\ref{eqn:CPW1}-\ref{eqn:CPWn}) are complete for any $\mu \geq \mu_0$.
\end{conjecture}

The proposed method for constructing CPWs in high dimensions is essentially different from the usual way of generalizing a 1D basis to a multidimensional basis by using the tensor product. 
Moreover, it is also clear that the index $n$ controls the size of the compact support and the scale of CPWs, while $\j$ controls the shift. These two parameters are analogues of the scale and shift parameters in  wavelet theory, which might help us build a new method of multiresolution analysis in our future work.
 In addition, the following scaling formula indicates the relation between the parameters $\mu$ and $\w$, which can be easily obtained by change of variables. 
 \begin{property}[Scaling formula] If we write $\psi^n_{\{\mu,\w\}}(\x)$ as the $n$-th BCPW obtained from Eqs. (\ref{eqn:CPW1}-\ref{eqn:CPWn}) with parameter $\mu,\w$, % (and with $V=0$), 
 then the following scaling formula holds,
 \begin{equation}
 \psi^n_{\{\mu,\w\}}(\x) = s^{d/2}\psi^n_{\{s^{2+d/2}\mu,s\w\}}(s\x).
 \end{equation}
 \end{property}
 
\subsection{Numerical algorithms}
\label{CPW:numerical}
To simplify our discussion, we only consider $\Omega=[0,L]$ in 1D with periodic boundary conditions. However, the algorithms for the fast transform and the fast inverse transform discussed below can be easily extended to higher dimensions. 

To solve the proposed the CPWs, we first solve Eq. (\ref{eqn:CPW1}). By introducing 
an auxiliary variable $u = \psi$, the constrained optimization problem [\ref{eqn:CPW1}] in 1D is equivalent to the following problem
\begin{eqnarray}
\psi^1 = \arg \min_{\psi,u} \frac{1}{\mu} \int |u| \d x  +  \int \psi \hat{H}_0\psi \d x \hspace{1.5cm} \nonumber \\
\text{s.t.} \quad  u = \psi \quad \&\quad \int \psi(x)\psi(x - j w)\d x = \delta_{j0}, \quad  j \in\mathbb{Z},
\end{eqnarray}
which can be solved by the following algorithm based on Bregman iteration (Refs.~\cite{Osher:2005,Yin:2008bregman,Goldstein:2009split}).

\begin{algorithm}
\label{alg:SOC_CPW}
Initialize $\psi^{1,0} = u^0 , b^0  = 0$.

\While{``not converged"}{
\begin{enumerate}
\item $\displaystyle \psi^{1,k} = \arg\min_{\psi} \int \psi \hat{H}_0 \psi \d x + \frac{\lambda}{2}\int (\psi - u^{k-1} + b^{k-1})^2\d x$ \, 
        \text{\rm s.t.} \, $\displaystyle \int \psi(x)\psi(x - j w)\d x = \delta_{j0}, \quad  j \in\mathbb{Z}$.
\item $
         \displaystyle u^k = \arg\min_{u} \frac{1}{\mu}\int  |u| \d x + \frac{\lambda}{2}\int (\psi^{1,k} - u + b^{k-1})^2 \d x
             = \text{\rm Shrink}(\psi^{1,k}+ b^{k-1},\frac{1}{\lambda\mu})$
\item $b^{k} = b^{k-1} +  \psi^{1,k} - u^k$.
\end{enumerate}
}
\end{algorithm}

\begin{figure}[h]
\centering
\begin{minipage}{0.49\linewidth}
\centering
\includegraphics[width=1\linewidth]{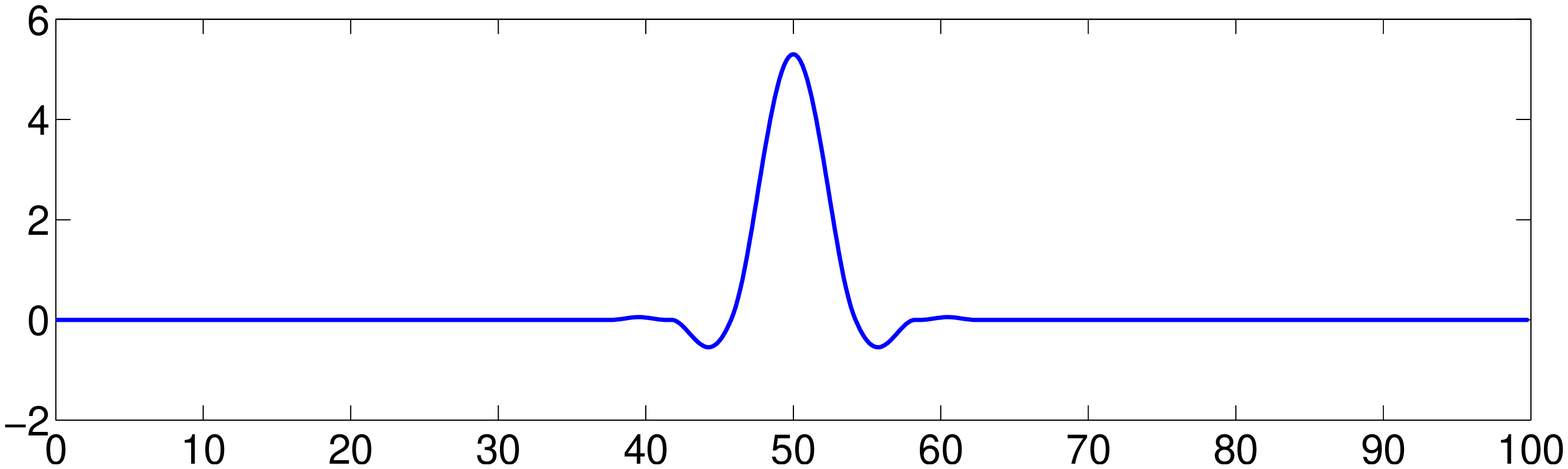}\\
$\psi_1$ \\
\end{minipage}\hfill
\begin{minipage}{0.49\linewidth}
\centering
\includegraphics[width=1\linewidth]{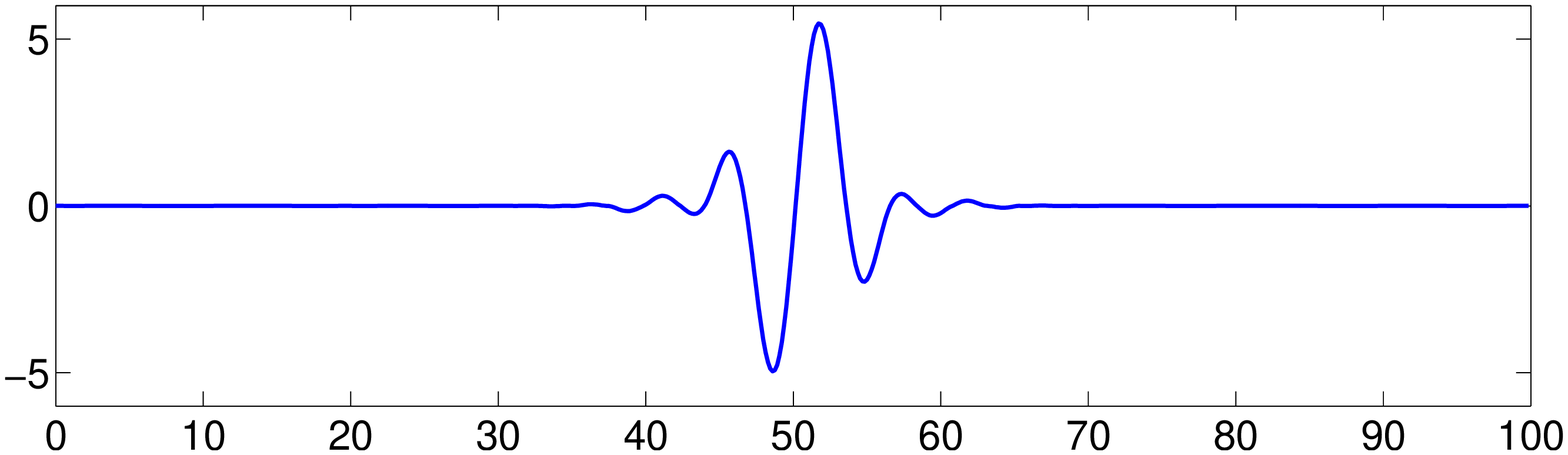}\\
$\psi_2$ \\
\end{minipage}\hfill\\
\begin{minipage}{0.49\linewidth}
\centering
\includegraphics[width=1\linewidth]{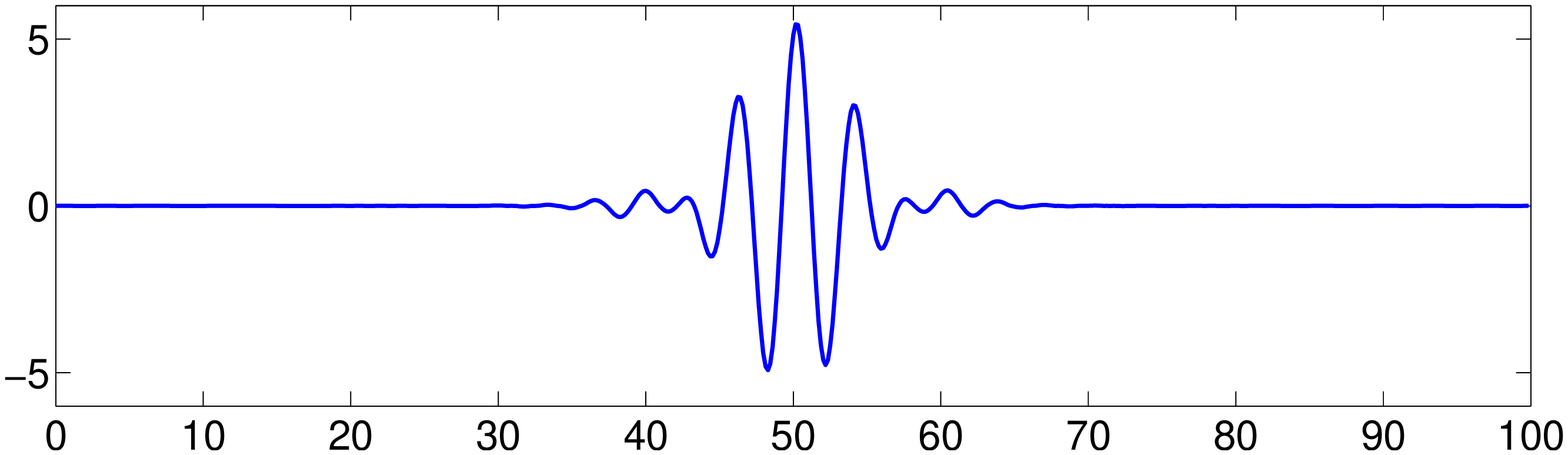}\\
$\psi_3$ \\
\end{minipage}\hfill
\begin{minipage}{0.49\linewidth}
\centering
\includegraphics[width=1\linewidth]{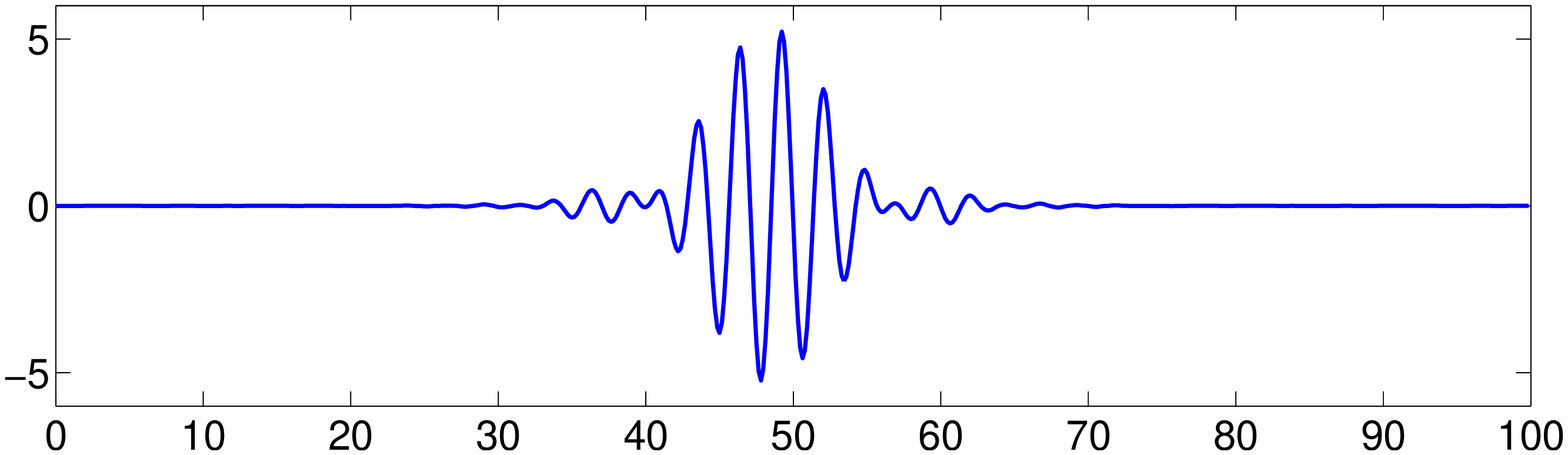}\\
$\psi_4$ 
\end{minipage}\hfill
\begin{minipage}{0.49\linewidth}
\centering
\includegraphics[width=1\linewidth]{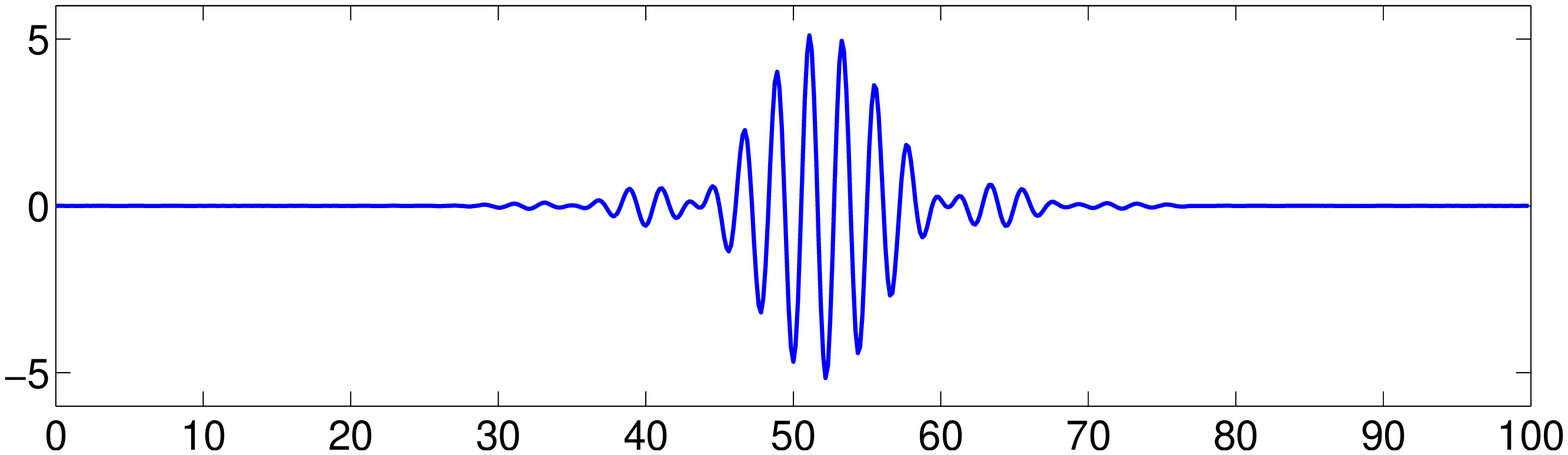}\\
$\psi_5$ 
\end{minipage}\hfill
\begin{minipage}{0.49\linewidth}
\centering
\includegraphics[width=1\linewidth]{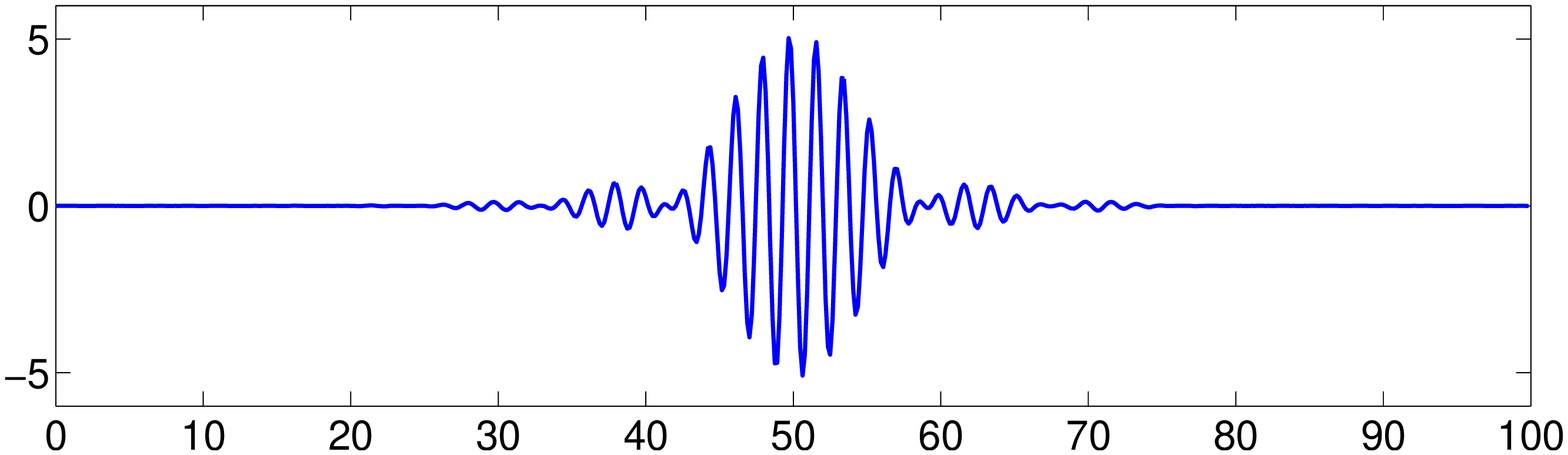}\\
$\psi_6$ 
\end{minipage}\hfill
\caption{From top to the bottom, the first six modes $\psi^1,\psi^2,\psi^3,\psi^4,\psi^5,\psi^6$ obtained by Eq. (\ref{eqn:CPW1} -\ref{eqn:CPWn}) using $L=100, \mu = 5,w = 5$.}
\label{fig:CPW}
\end{figure}

We solve the first problem in the above algorithm in the Fourier space, since the kinetic energy and the constraints are diagonal. In other words, let's write $\displaystyle \psi = \sum_{G} \psi_G e^{iG x}, u^{k-1} = \sum_{G} u^{k-1}_G e^{iG x}$ and $b^{k-1} = \sum_{G} b^{k-1}_G e^{iG x}$. Then we need to solve the following problem
\begin{eqnarray}
 \min_{\{\psi_G\},\gamma_j}\sum_{G}\frac{G^2}{2}\psi_G^2 
 + \frac{\lambda}{2}\sum_{G} (\psi_G - u_G^{k-1} + b_G^{k-1})^2    
 + \sum_{j} \gamma_j (\sum_{G}\cos(G j w)\psi_G^2 - \delta_{j0}),      
\label{eqn:constraint_fs}
\end{eqnarray}
where $\{\gamma_j\}$ are Lagrangian multipliers associated with the orthonormality constraints. Then $\{\gamma_j\}$  are found by solving the following equations:
\begin{equation}
\left\{\begin{array}{c}
\displaystyle\sum_{G}\cos(G j w)\psi_G^2  = \delta_{j0}, \quad\quad  j \in\mathbb{Z} \\
\displaystyle\psi_{G} = \frac{\lambda(u_G - b_G) }{\lambda + G^2 + 2\sum_j\cos(G j w) \gamma_j}
\end{array}\right. 
\end{equation}
One can go further and define higher order modes $\psi^n$  that satisfy Eqs. (\ref{eqn:CPW1}-\ref{eqn:CPWn}). The additional orthogonality constraints can be handled via one more orthogonality constraints splitting using a similar way as in Algorithm \ref{alg:CM_SOC}  in the manner of SOC proposed in Ref.~\cite{Lai:2013JSC}. 
In Figure~\ref{fig:CPW}, we illustrate the first six modes $\psi^1,\psi^2,\psi^3,\psi^4, \psi^5, \psi^6$ obtained from Eqs.~(\ref{eqn:CPW1}-\ref{eqn:CPWn}) using $L=100, \mu = 5,w = 5$.

\begin{figure}[htp]
\centering
\begin{minipage}{0.49\linewidth}
\includegraphics[width=1\linewidth]{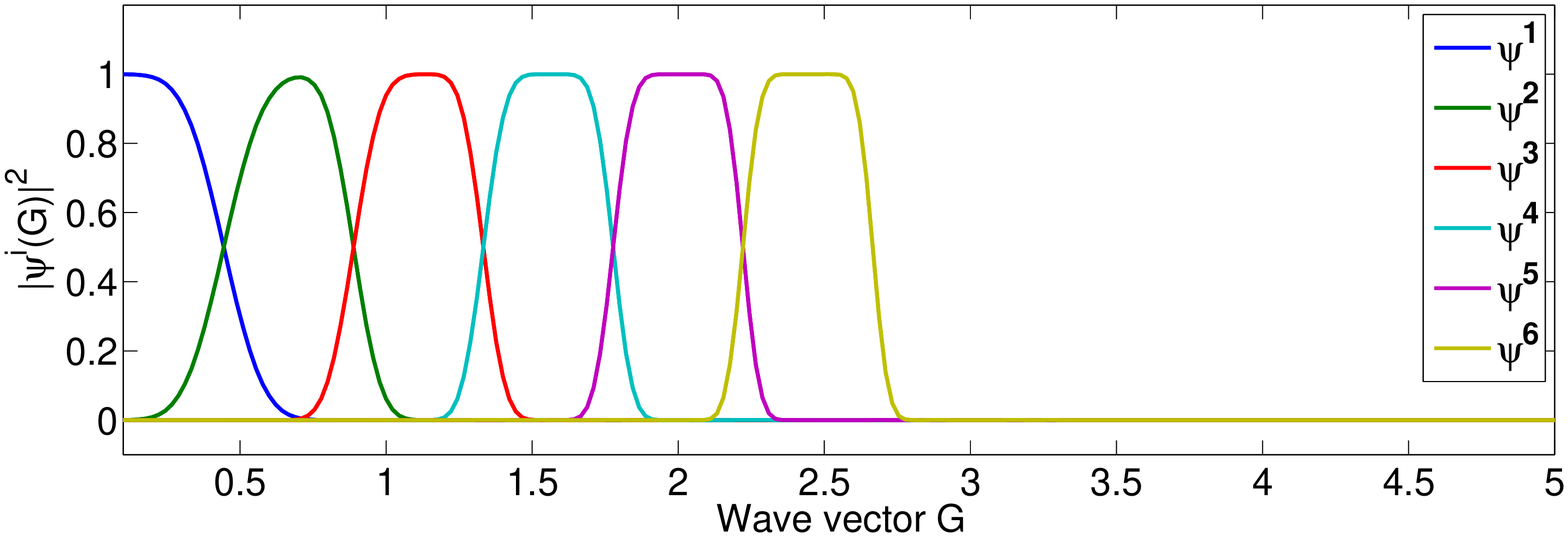}
\end{minipage}\hfill
\begin{minipage}{0.49\linewidth}
\includegraphics[width=1\linewidth]{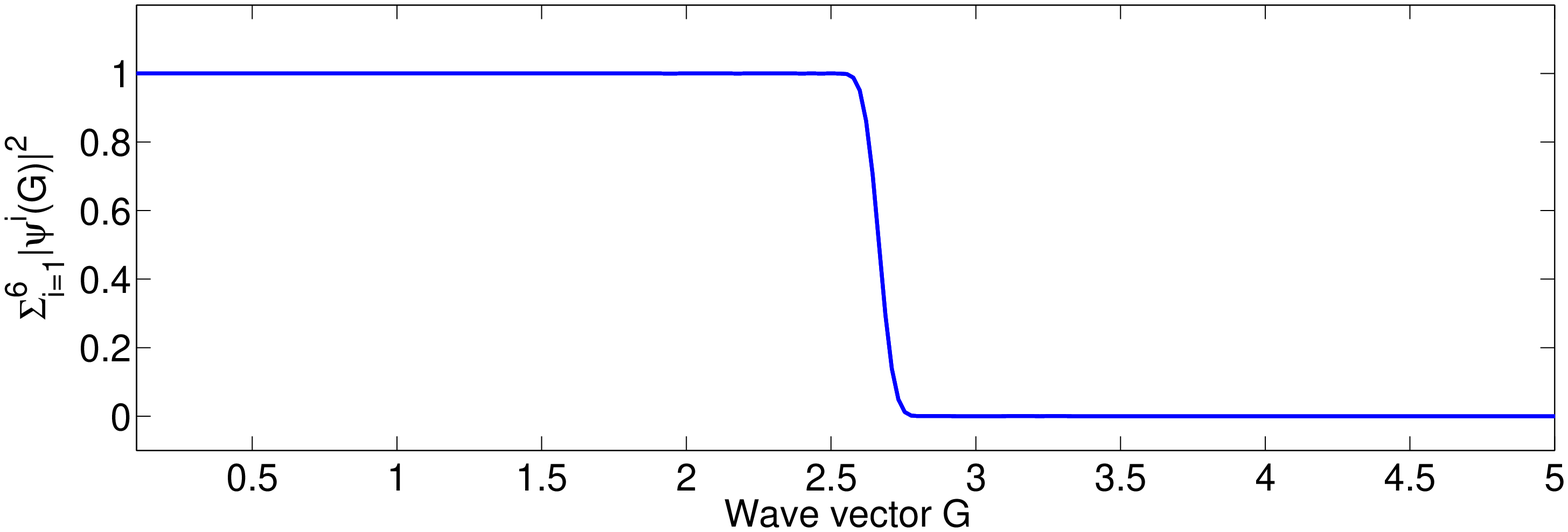}
\end{minipage}\hfill
\caption{Spectral density distribution of CPWs. Left: the spectral density distribution of $\psi^1,\psi^2,\psi^3,\psi^4,\psi^5,\psi^6$. Bottom:  The total spectral density distribution of the first four modes.}
\label{fig:CPW_spectraldensity}
\end{figure}

The next interesting property is the distribution of the spectral weight $|\psi^i_G|^2$ and the total spectral weight $\sum_{i=1}^N|\psi^i_G|^2$ of the first six BCPWs. 
In Figure~\ref{fig:CPW_spectraldensity}, we show the spectral weight distribution of the first six BCPWs and their total spectral weight. From the left panel of Figure~\ref{fig:CPW_spectraldensity}, we can see that each mode occupies a distinct region in the Fourier space. Moreover, the total spectral weight of the first six BCPWs forms a nice-looking step function, which is a desirable property for obtaining convergence rates similar to the plane wave basis. In other words, locally the basis covers approximately the same Fourier space as plane waves below a given kinetic energy cutoff. 

\subsection{CPW representations for localized functions}
\label{CPW:rep}

\begin{figure}[ht]
\centering
%\begin{minipage}{0.49\linewidth}
\includegraphics[width=.8\linewidth]{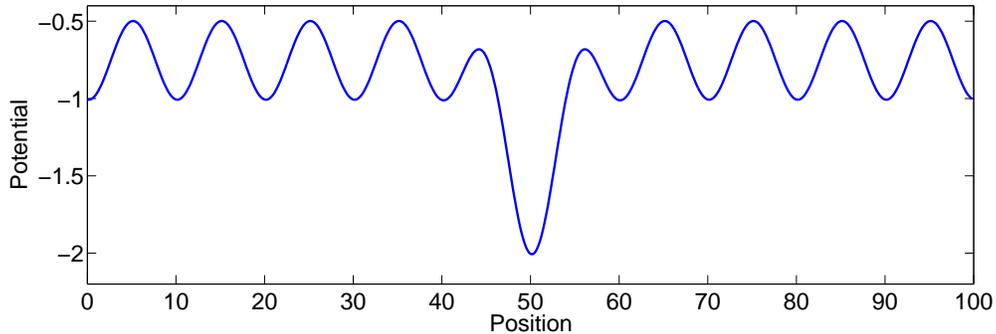}\\
%\end{minipage}\hfill
\caption{The potential function $V(x)$ of the impurity Kronig-Penney model.}
\label{fig:impurityKP}
\end{figure}

Localized properties of the proposed CPWs can be expected to have advantages in describing local fine structures in many applications. 
For example, localized functions can be obtained from energy state functions of  a physical system with one of the potential wells deeper than the others. This is analogous to introducing an  ``impurity" atom in quantum mechanics as the potential function illustrated in Figure~\ref{fig:impurityKP}, which we call an ``impurity'' Kronig-Penny (IKP) model. Using the first 120 CPWs generated by the first six BCPWs illustrated in Figure~\ref{fig:CPW}, we would like to demonstrate several advantages of the CPW representation for localized functions.

\begin{figure}[h]
\centering
\begin{minipage}{0.49\linewidth}
\centering
\includegraphics[width=1\linewidth]{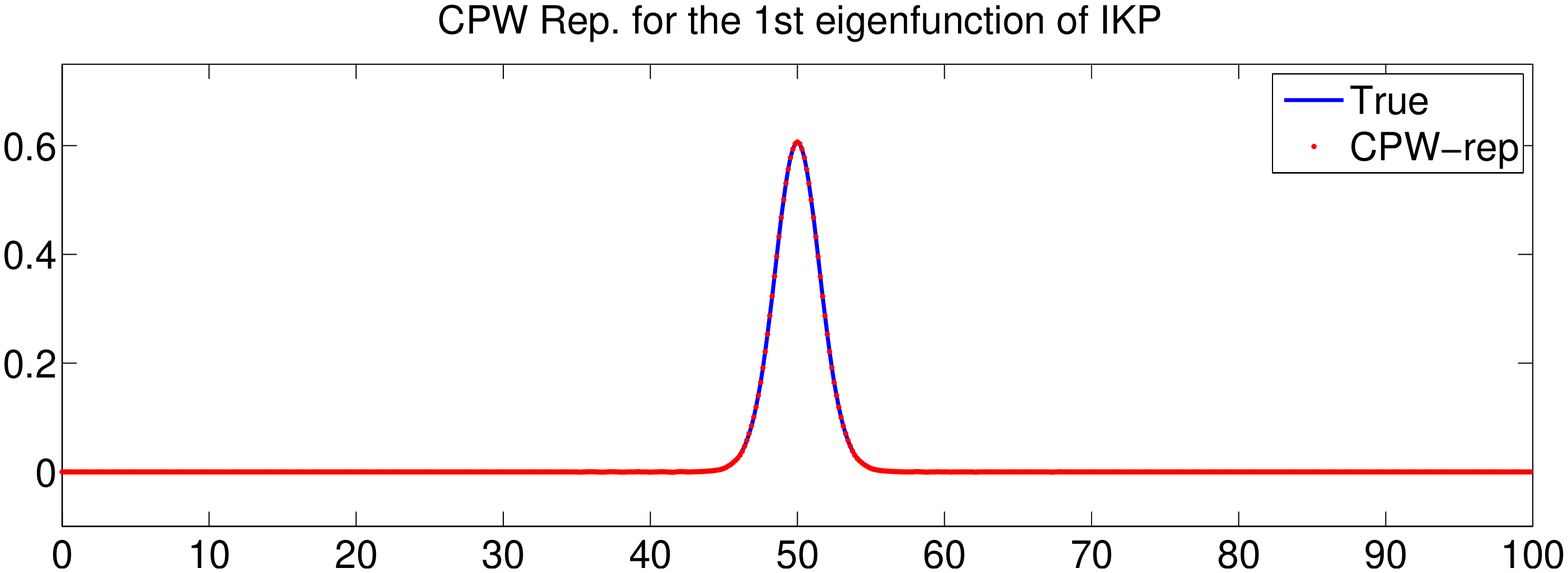}\\
\vspace{0.2cm}
\includegraphics[width=1\linewidth]{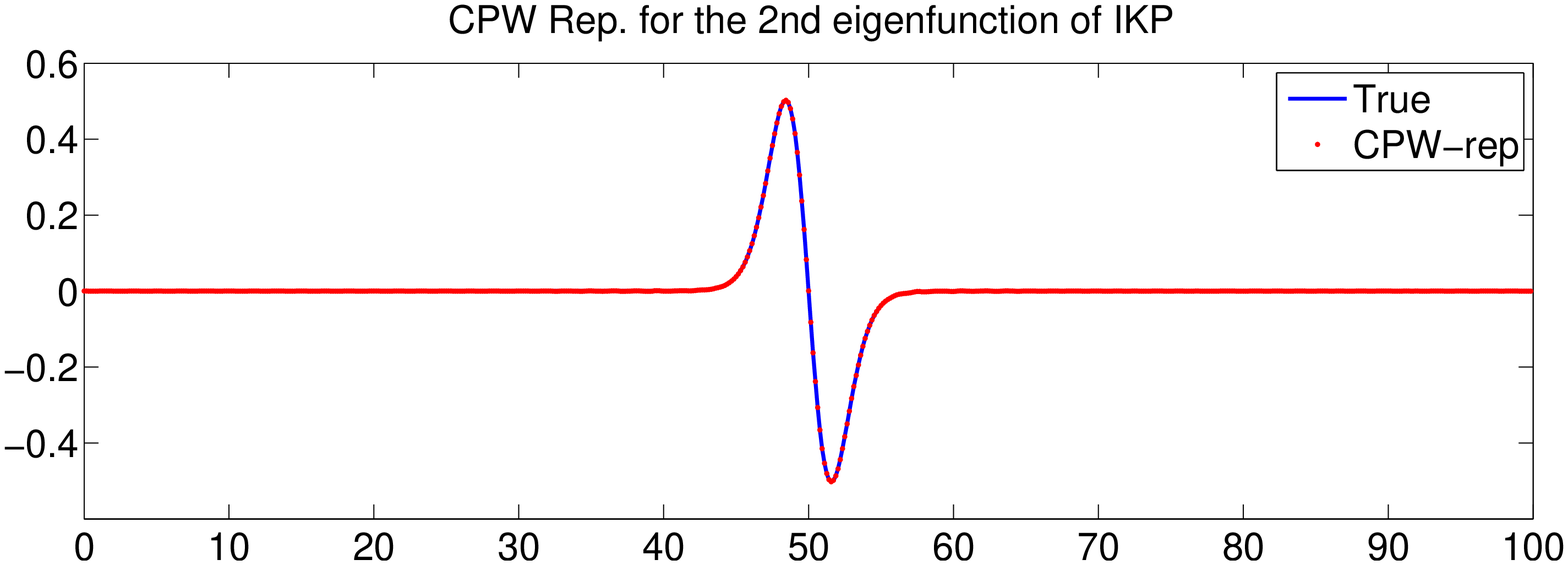}\\
\vspace{0.2cm}
\includegraphics[width=1\linewidth]{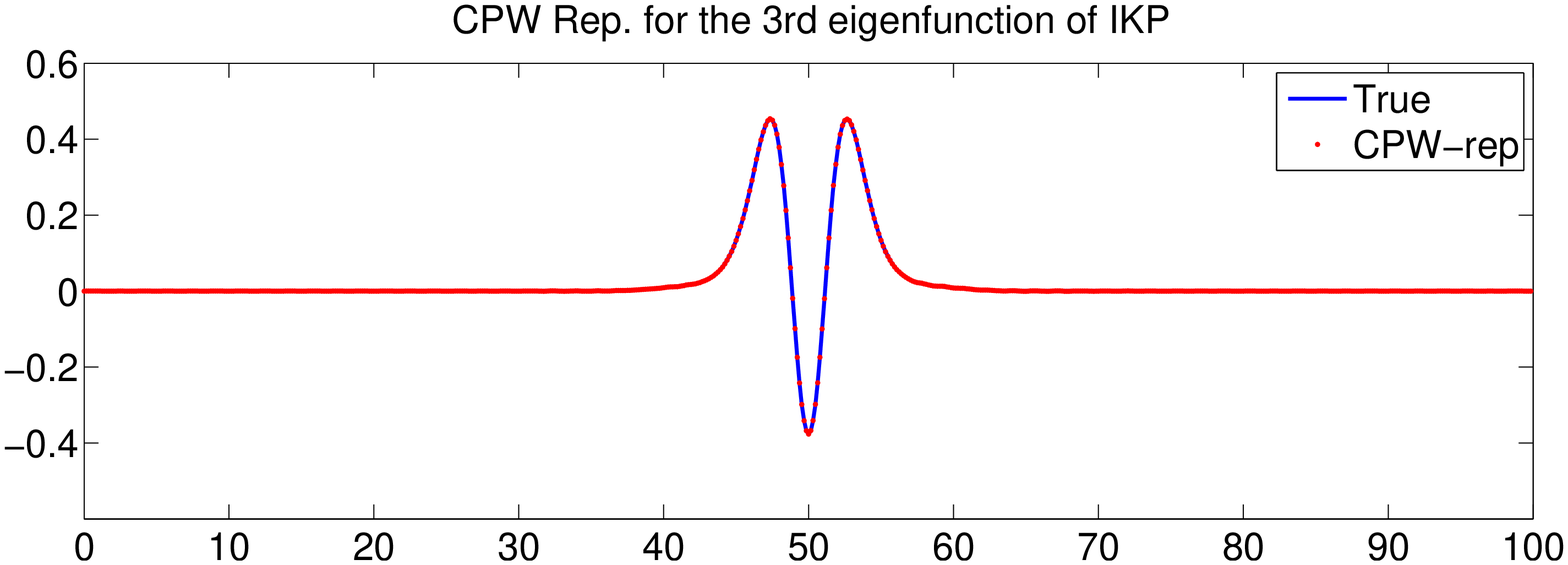}\\
\vspace{0.2cm}
\includegraphics[width=1\linewidth]{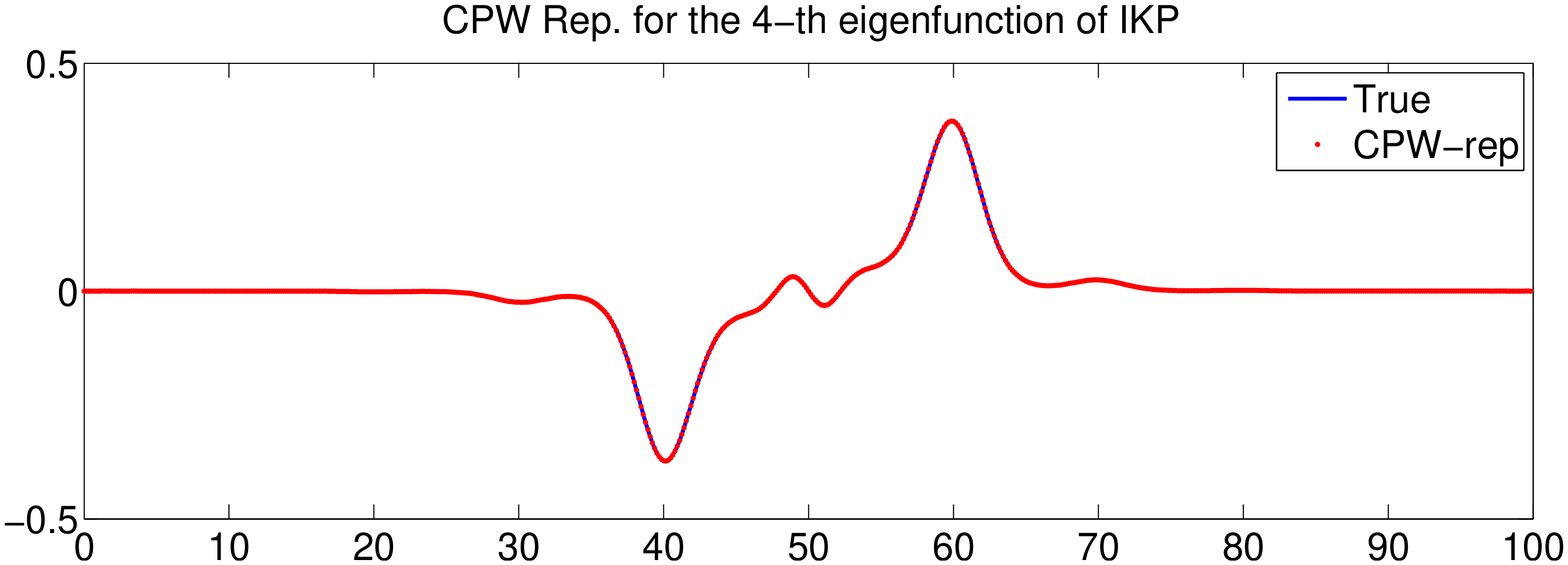}\\
\end{minipage}\hfill
\begin{minipage}{0.49\linewidth}
\centering
\includegraphics[width=1\linewidth]{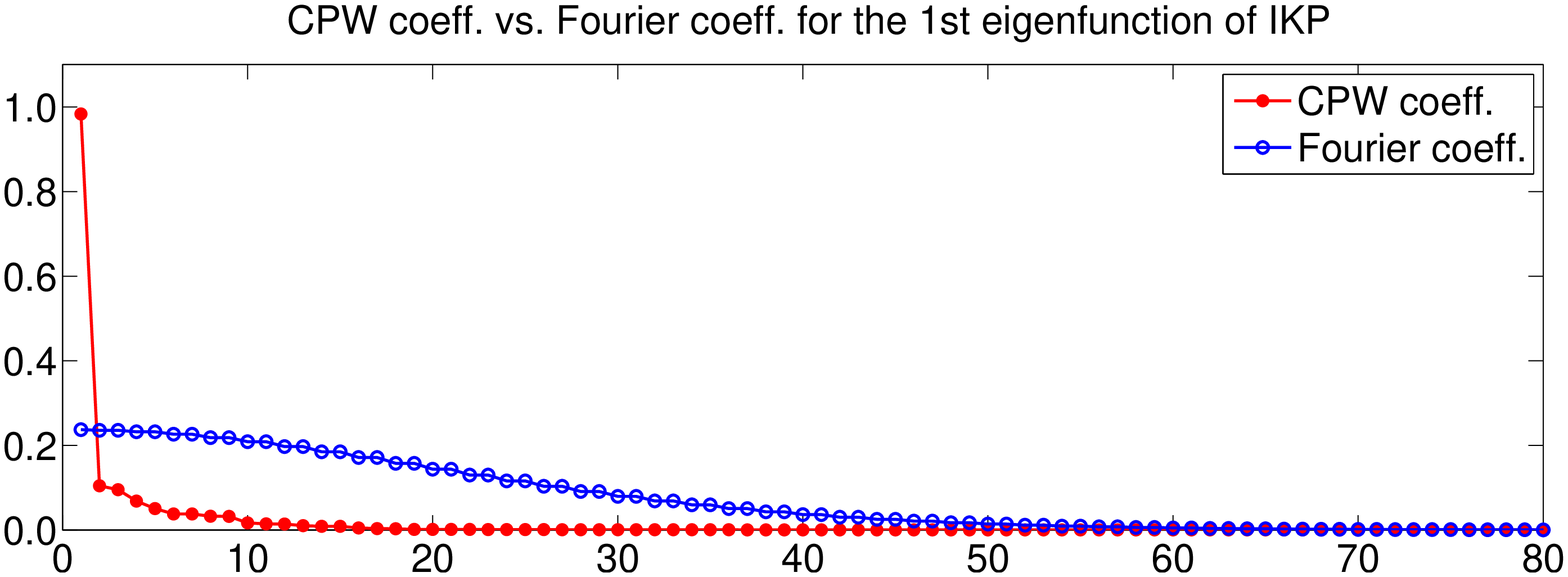}\\
\vspace{0.2cm}
\includegraphics[width=1\linewidth]{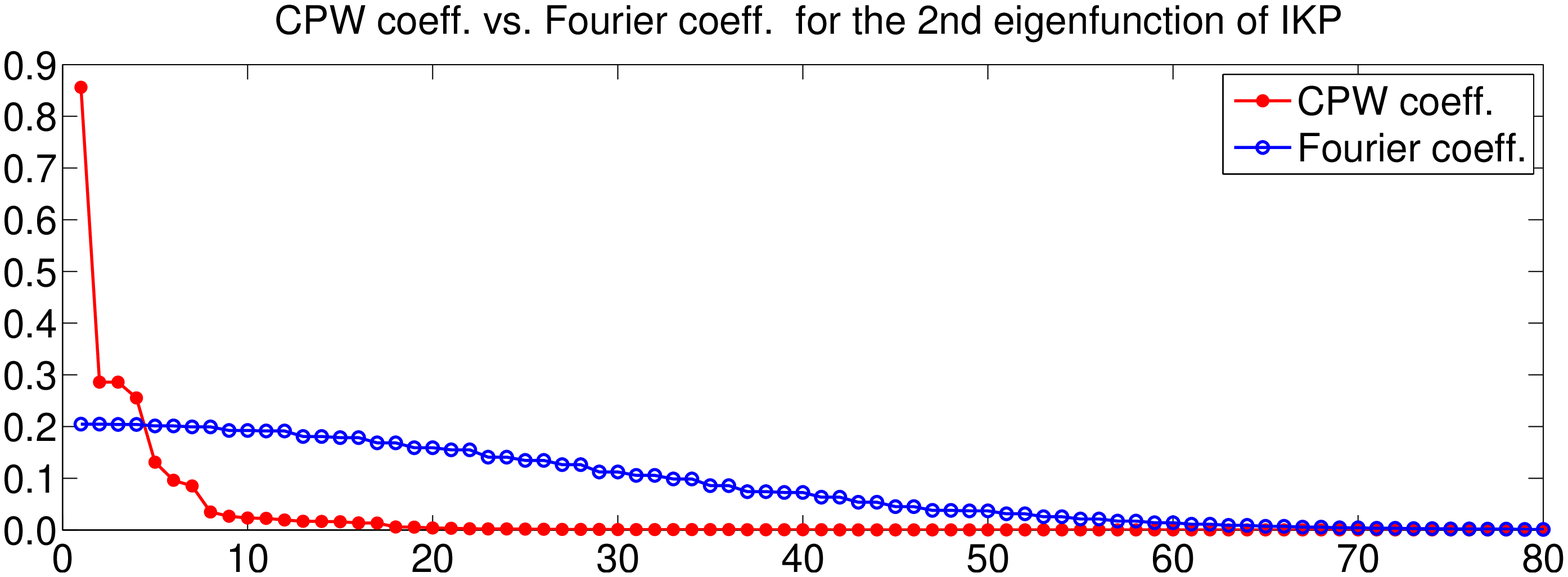}\\
\vspace{0.2cm}
\includegraphics[width=1\linewidth]{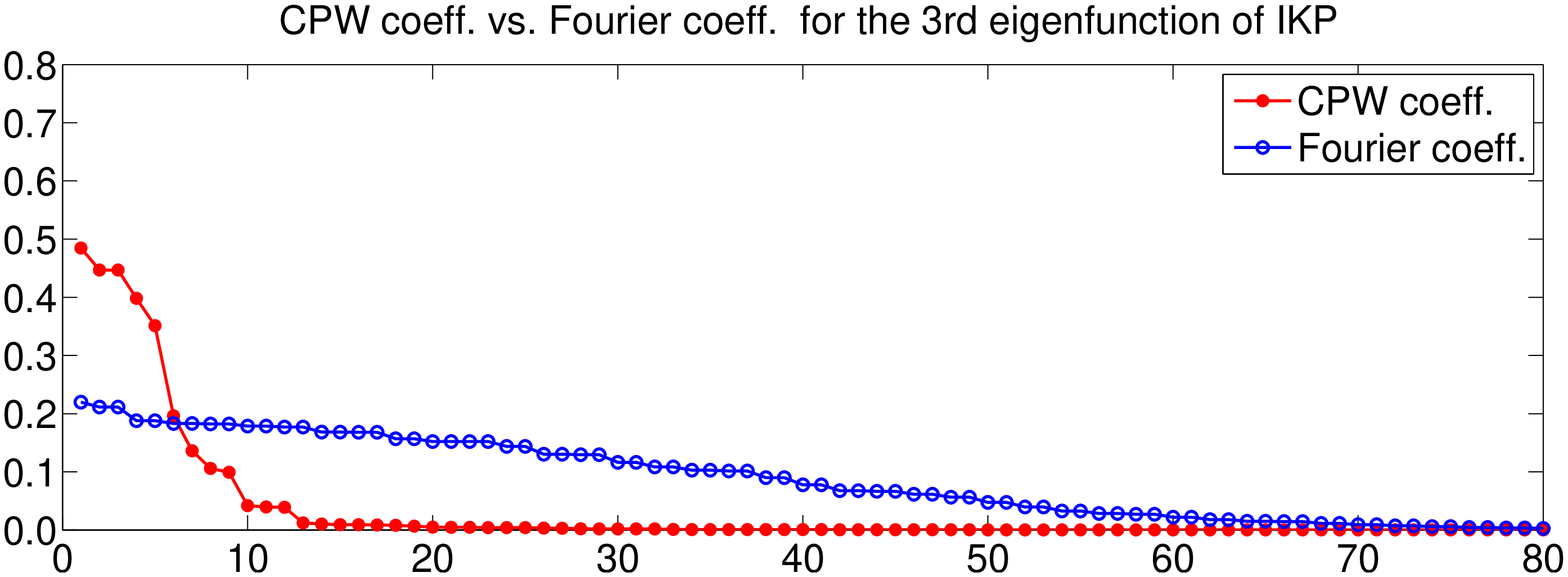}\\
\vspace{0.2cm}
\includegraphics[width=1\linewidth]{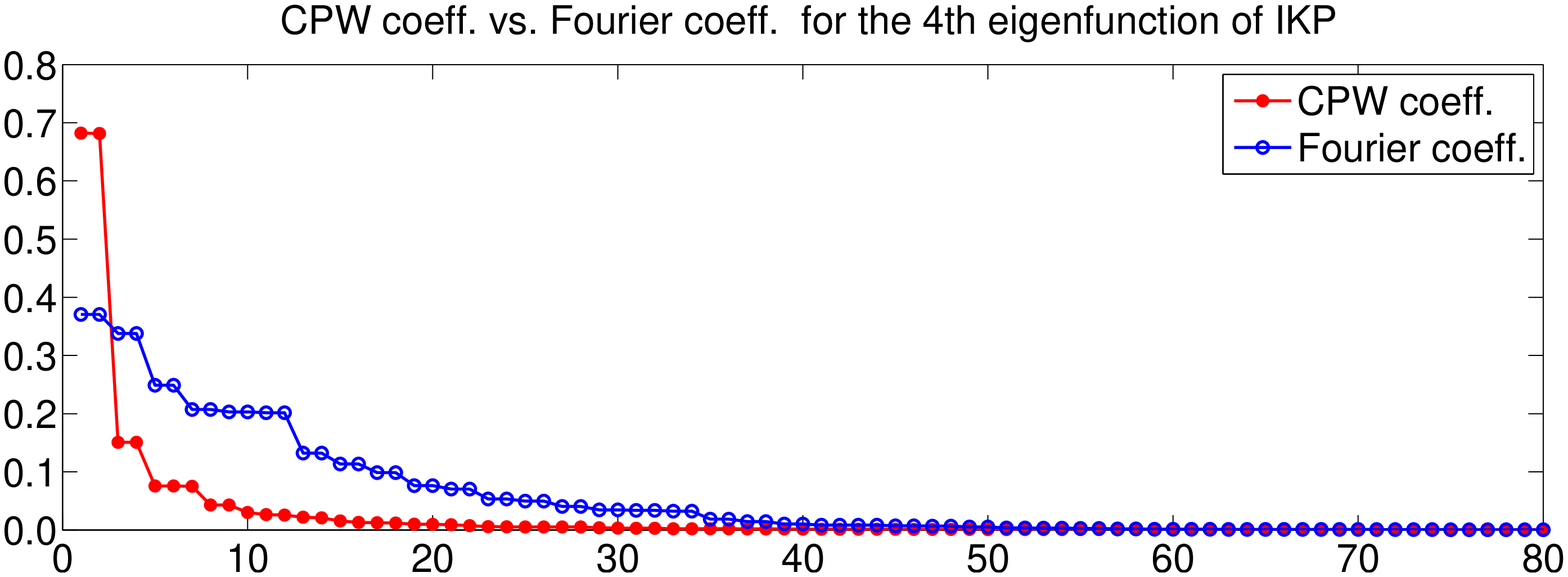}\\
\end{minipage}\hfill
\caption{CPW representations for the first four energy states of the IKP model. Left: Comparisons of CPW representations for the first four eigenfunctions of the IKP model with their true values. Right: Comparisons of the first 80 largest magnitude coefficients of CPW representations and classical Fourier function representations.}
\label{fig:CPW_rep}
\end{figure}

Several localized energy states can be obtained from the IKP model. The blue curves in the left four images of Figure~\ref{fig:CPW_rep} illustrate the first four energy states of the IKP model. We can successfully recover these localized functions using the first 120 CPWs generated in Section \ref{CPW:numerical}. Representation results plotted by red dots in the left four images of 
Figure~\ref{fig:CPW_rep} demonstrate the accuracy of CPW representations. In addition, we also compare the first 80 largest magnitude coefficients of CPW representations and classical Fourier function representations for these four functions. The right column in Figure~\ref{fig:CPW_rep} clearly shows that CPWs can provide a much more sparser representation for these four localized functions than the Fourier basis. 
%
% we compare the first 120 CPWs representations for the first six eigenfunctions of the IKP model with their true functions. Our experiments show that only 25 CPWs can accurately recover the first few eigenfunctions of  the IKP model, while 73 Fourier basis functions are needed for the same accuracy. 
Table~\ref{tab:IKP_RepError_CPW_Fourier} also reports more detailed comparisons of the  $\ell_2$ error of the CPW representation and Fourier function representation using the first few largest magnitude coefficients. 
It is promising that the number of CPWs needed to accurately represent the localized functions is significantly smaller than the number of Fourier basis functions for the same accuracy.

%\begin{figure}[t]
%\centering
%\begin{minipage}{0.49\linewidth}
%\includegraphics[width=1\linewidth]{IKP_CPW_Rep_error.eps}\\
%\end{minipage}\hfill
%\begin{minipage}{0.49\linewidth}
%\includegraphics[width=1\linewidth]{IKP_Fourier_Rep_error.eps}\\
%\end{minipage}\hfill\\
%\centering
%\resizebox{6.5in}{1.5cm}{
\begin{table}
\centering
\begin{tabular}{|c||c|c|c|c||c|c|c|c|c|}
\hline
\multirow{2}{1.1cm}{$\#$ of modes} &  \multicolumn{4}{|c||}{Rep. Error using CPWs} & \multicolumn{4}{|c|}{Rep. Error using Fourier basis}   \\
 \cline{2-9}
 &  Error$_{f_1}$ & Error$_{f_2}$ &Error$_{f_3}$ &Error$_{f_4}$  & Error$_{f_1}$ & Error$_{f_2}$ &Error$_{f_3}$ &Error$_{f_4}$ \\
 \hline
20 &   0.0091   & 0.0170  &  0.0305 &   0.0506   & 1.0282 &   1.3586  & 1.4886  &  0.4553   \\
 \hline
30   & 0.0051   &  0.0084  &  0.0111  &  0.0204    & 0.5052 &    0.8289  &  0.9951  &  0.2004  \\
 \hline
40   & 0.0038  &  0.0065  &  0.0076  &  0.0115    & 0.2112  &  0.4332 &   0.5875  &  0.0619   \\
     \hline
50   & 0.0036  &  0.0063 &   0.0066  &   0.0075   & 0.0763  &  0.1930  &  0.2964  &  0.0228 \\
    \hline
60   & 0.0035  &  0.0063  &  0.0064   &  0.0052  & 0.0242 &   0.0651  &  0.1304 &   0.0097 \\
    \hline
70   & 0.0035  &  0.0063  &  0.0064   &  0.0043    & 0.0068 &   0.0192  &  0.0470  &  0.0048 \\
    \hline
\end{tabular}
\caption{Comparisons of Representation error using CPWs and Fourier basis.}
\label{tab:IKP_RepError_CPW_Fourier}
\end{table}
%\caption{Representation error comparisons using CPWs and Fourier basis. Top figures: representation error curves of the first six eigenfunctions %of the IKP model using CPWs (left) and Fourier basis functions (right) via different numbers of modes. Bottom table: Representation errors using %CPWs vs. Fourier basis.}
%\label{fig:IKP_RepError_CPW_Fourier}
%\end{figure}

\begin{figure}[ht]
\centering
\begin{minipage}{0.49\linewidth}
\includegraphics[width=1\linewidth]{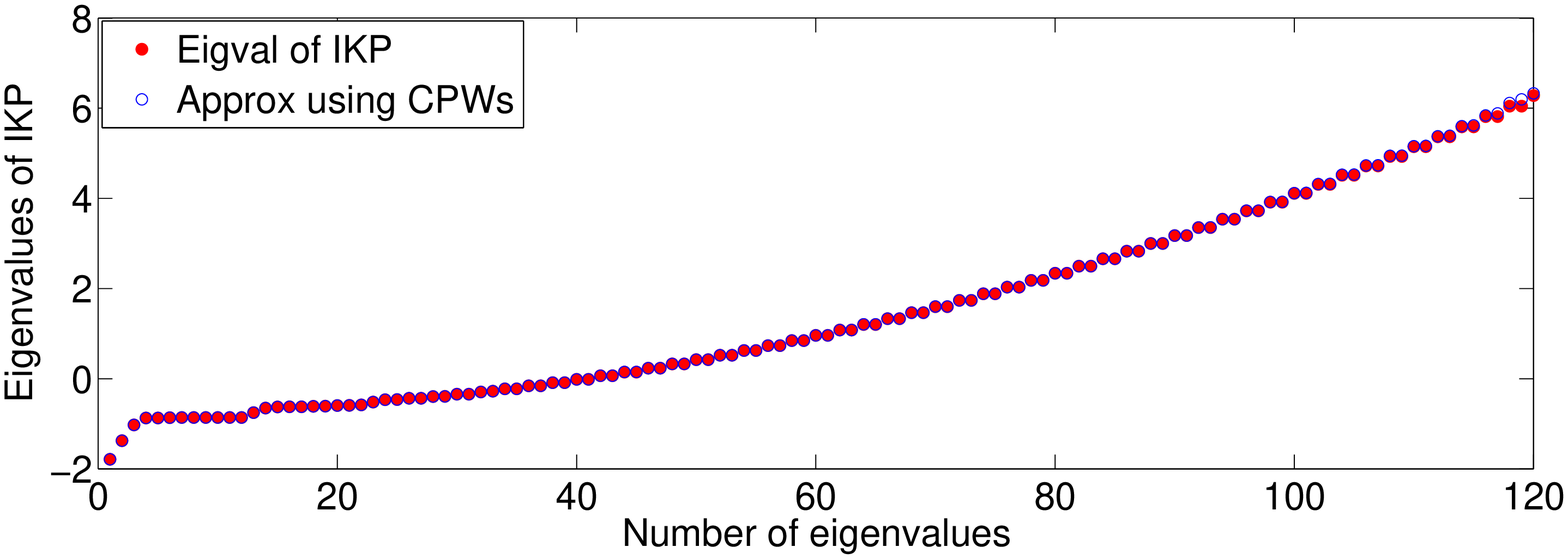}\\
\end{minipage}\hfill
\begin{minipage}{0.49\linewidth}
\includegraphics[width=1\linewidth]{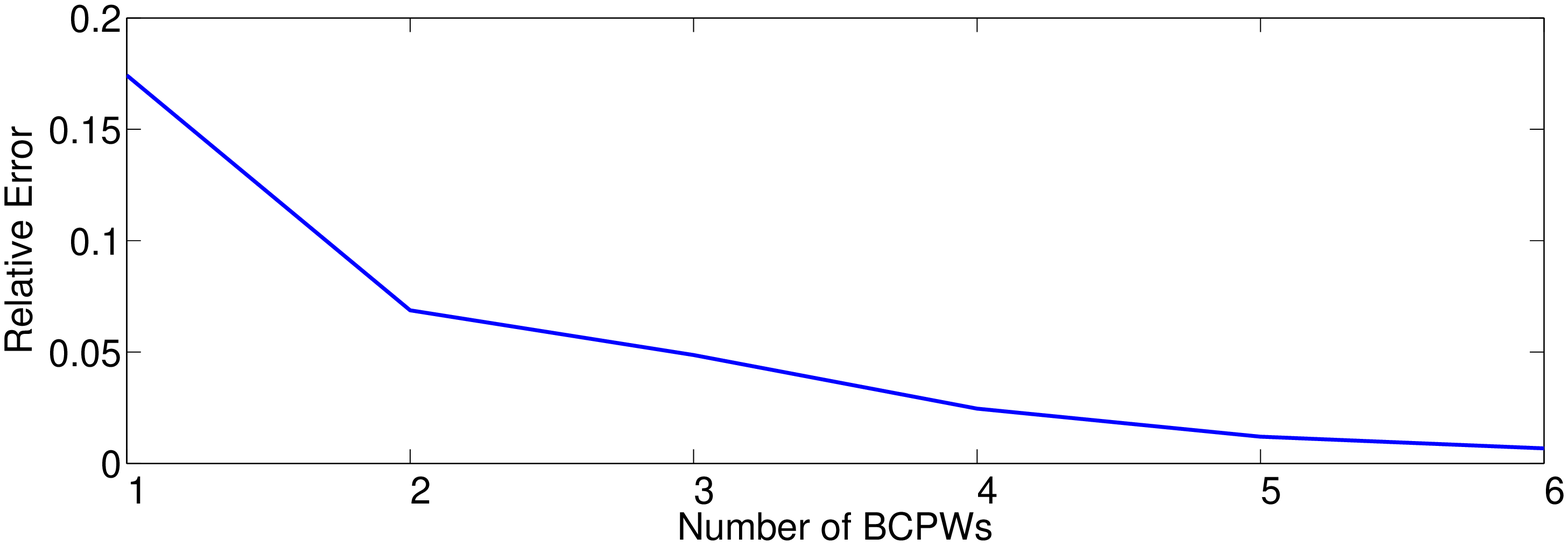}\\
\end{minipage}\hfill
\caption{Eigenvalues approximation for the IKP model using the CPW representation. Left: the first 120 approximated eigenvalues vs. true values. Right: relative eigenvalue approximation error using CPWs generated by the first 1 to 6 BCPWs. }
\label{fig:IKP_eigs}
\end{figure}

Due to the advantage of representing localized functions using CPWs,  we can also successfully approximate the first few eigenvalues of the IKP model using CPWs. In other words, we use eigenvalues of the matrix $\left(\langle \b_j^n, \hat{H} \b_i^m \rangle\right)_{i,j=1,\cdots 20}^{n,m = 1,\cdots  N_{BCPWs}}$ ($N_{BCPWs}$ is the number of BPCWs) to approximate the first few eigenvalues of the Schr\"odinger operator $\hat{H}$ used in IKP model. 
In the left image of Figure~\ref{fig:IKP_eigs}, the red dots plot the first 120 eigenvalues of IKP model obtained by using spectral method with 640 nodes, while the blue circles are approximation results using the first 120 CPWs generated by the first 6 BCPWs. The right error curve in Figure~\ref{fig:IKP_eigs} also reports the relative approximation error using CPWs generated by the first 1 to 6 BCPWs. It is clear that CPW representations provide accurate approximations for the first few eigenvalues of $\hat{H}$, where original eigenvalue problem can be reduced to an eigenvalue problem of a small size matrix. 

The conclusion is that CPWs become attractive if we need to represent a function that varies slowly through most of the space, except for a few regions; plane waves would need to increase the energy resolution everywhere in space, while we can do that locally with CPWs. 

Moreover, representing functions using the CPWs can be efficiently processed. In the rest of this section, we introduce a fast CPW transform and a fast inverse CPW transform.

\subsection{Fast CPW transform}

Given a function $f \in L^2([0,L])$, we have a transformation from $f$ in real space to the basis coefficients $\{f^n_j\}$ in frequency space. Recall that:
\begin{equation}
\displaystyle f^n_j = \int f(x) b^n_j(x)  = \sum_G (b^n_G)^* e^{-i(G j w)}\int f(x) e^{-i Gx}.
\end{equation}
where $\{b^n_G\}$ are the coefficients of $b^n_0 = \psi^n$ in the Fourier space. In other words, we write $b^n_0 = \sum_{G} b^n_G e^{i G x}$, and we have
\begin{equation}
\displaystyle  b^n_j = \sum_{G} b^n_G e^{i G x} e^{i(G  j w)}.
\end{equation}
The computation of this transform can be performed efficiently in three steps as follows.
\begin{algorithm}[Fast CPW Transform]
\begin{enumerate}
\item Fourier Transform. $\displaystyle f_G=  \int f(x) e^{-iG x} \d x$.
\item Multiplication and summation: 
         $\zeta^n_m = \sum_{k} (b^n_{G_{m+kN_0}})^* f_{G_{m+kN_0}}, \quad m = 0,\cdots, N_0 -1$, where $N_0 = [L/\omega], G_m= \frac{2\pi m}{L}$.
\item Fourier transforms of length $N_0$ to get the basis coefficients: $\displaystyle f^n_j = \sum_{m=0}^{N_0 - 1} \zeta^n_m e^{-i G j w}$.
\end{enumerate}
\end{algorithm}

\begin{figure}[h]
\begin{minipage}{0.49\linewidth}
\centering
\includegraphics[width=1\linewidth]{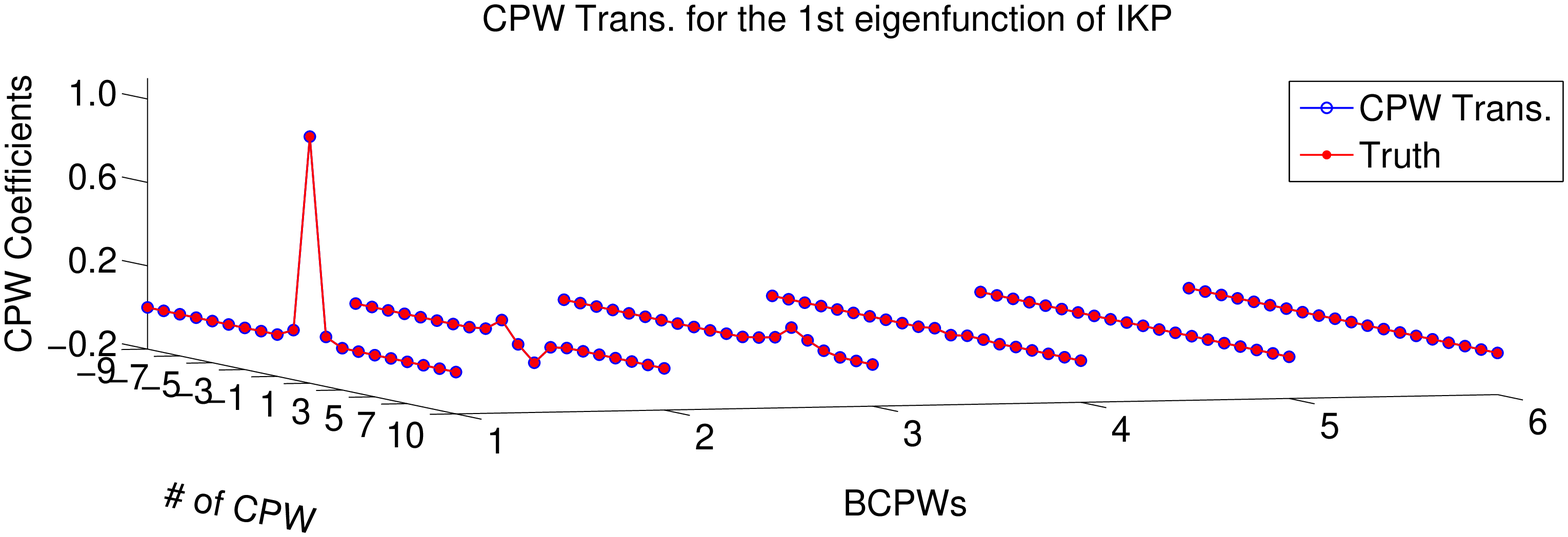}\\
\includegraphics[width=1\linewidth]{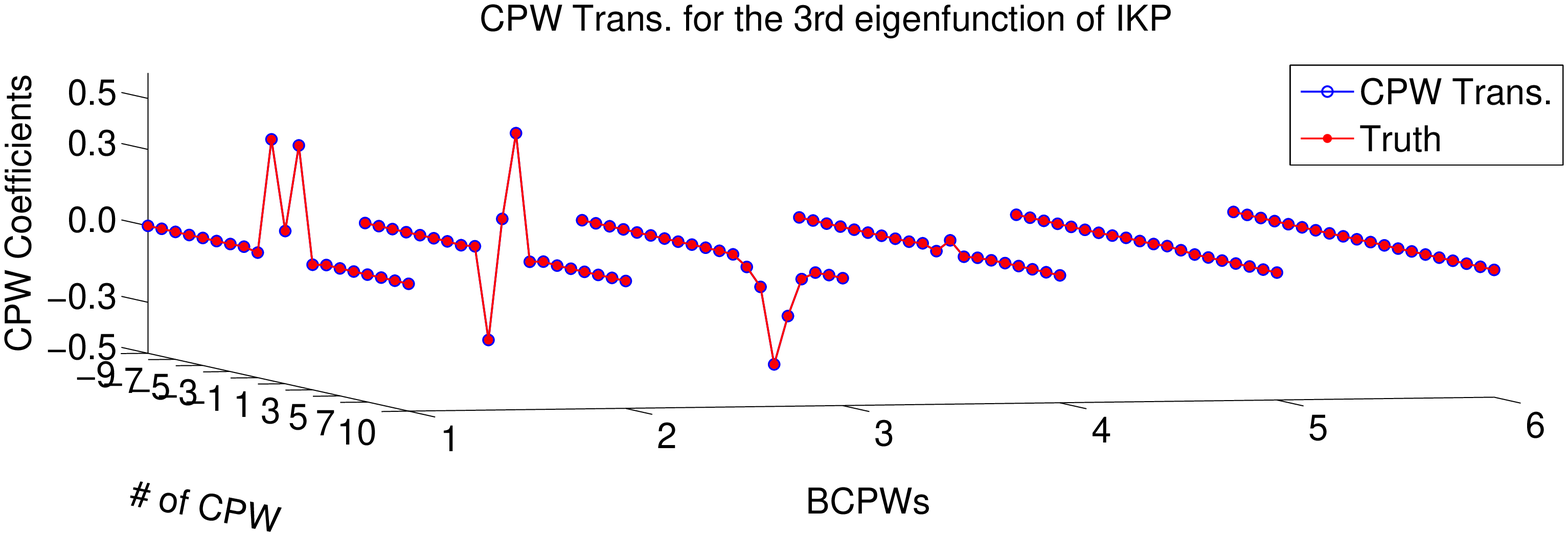}\\
\end{minipage}\hfill
\begin{minipage}{0.49\linewidth}
\centering
\includegraphics[width=1\linewidth]{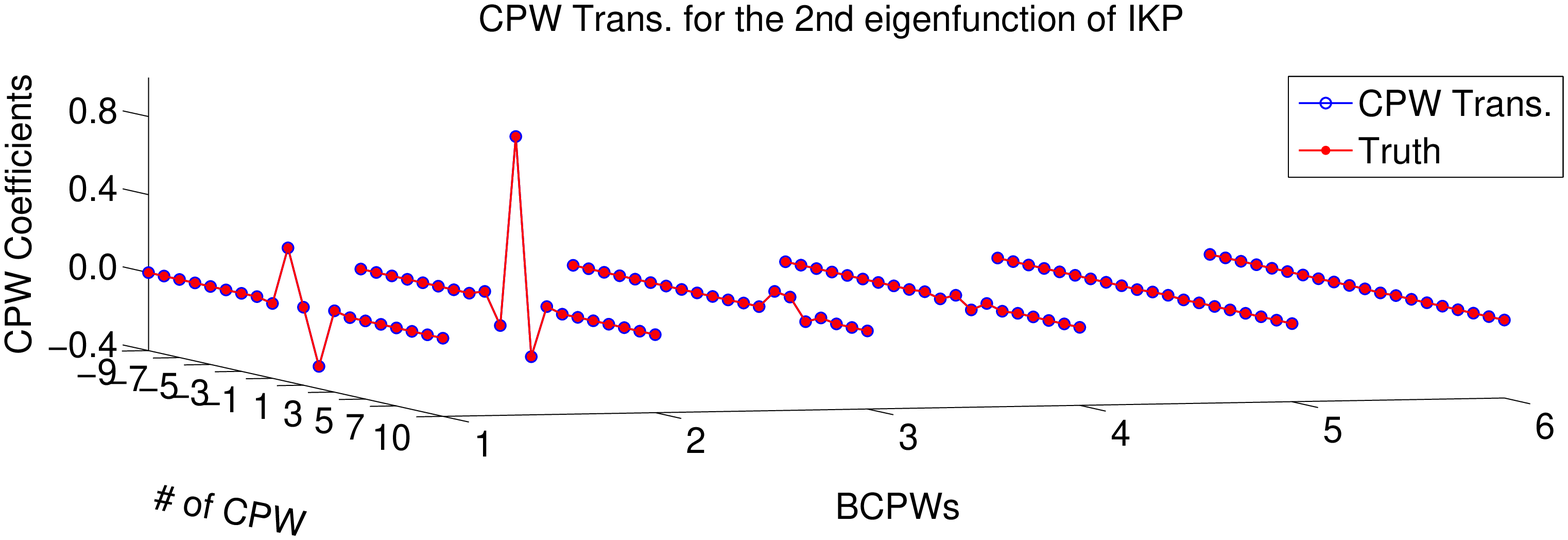}\\
\includegraphics[width=1\linewidth]{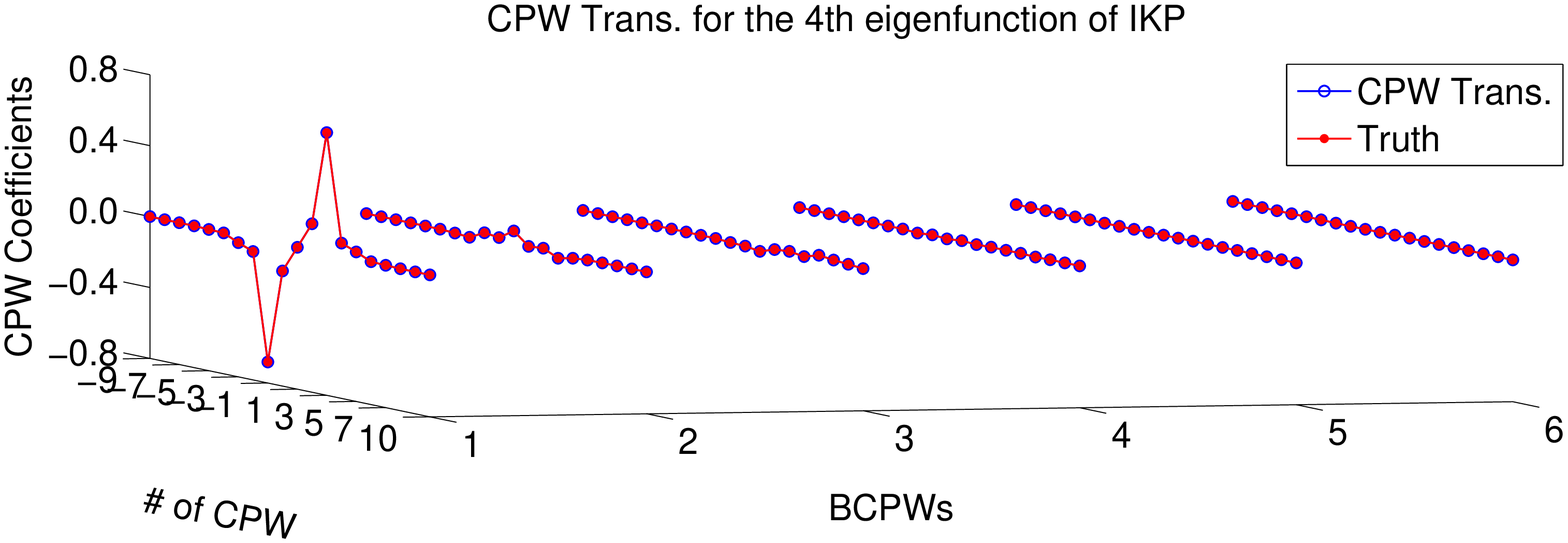}\\
\end{minipage}\hfill
\caption{CPW transform for  the first four lowest-energy states of the IKP model, where 6 levels of CPW transform are shown separately for each energy state. }
\label{fig:CPWtransform}
\end{figure}
We demonstrate this CPW transform for the first four lowest-energy states of the IKP model. Figure~\ref{fig:CPWtransform} illustrates accurate results using the CPW transform: the two sets of points are on top of each other, the red dots from direct diagonalization and the blue circles from the CPW transform described above.

Next, we can also have an inverse transformation from given CPW coefficients $\{f^n_j\}$ in frequency space to a function $f$ in real space. Recall that:
\begin{equation}
f(x) = \sum_{n} \sum_{j} f^n_j b^n_j(x) =  \sum_{n} \sum_{j} f^n_j \psi^n(x - j w).
\end{equation}

Therefore, one can rewrite $f$ as
\begin{equation}
f(x) =  \sum_{n,j} f^n_j b^n_j(x) =  \sum_{n} \sum_{G} b^n_Ge^{i G x} \sum_{j} f^n_j e^{i(G j w)}.
\end{equation}

The above summation can be efficiently computed in the following three steps.

\begin{algorithm}[Fast Inverse CPW Transform]
\label{alg:FastInverseTransform}
\begin{enumerate}
\item Fourier Transform. $\displaystyle \tilde{f}^n_m =  \sum_{j = 0}^{N_0 - 1} f^n_j e^{i(G_m j w)}$, where $N_0 = [L/\omega], G_m= \frac{2\pi m}{L}$.
\item Multiplication and summation. Note that all $\tilde{f}^n_m$ are periodic with a period $N_0$. Hence, we calculate
         $$  \xi_m^n = b^n_{G_m} \tilde{f}^n_{\text{\rm mod}(m,N_0)}.$$ Here, we only need to go up to Fourier coefficient values $m$ for which the corresponding basis function $b^n_{G_m}$ has nonzero coefficients. Then we add contributions from all $n$: $\displaystyle \xi_m = \sum_n \xi^n_m$.
\item Fourier transform to real space $\displaystyle f(x) = \sum_m \xi_m e^{i G_m x}$.
\end{enumerate}
\end{algorithm}

\begin{figure}[H]
\begin{minipage}{0.49\linewidth}
\centering
\includegraphics[width=1\linewidth]{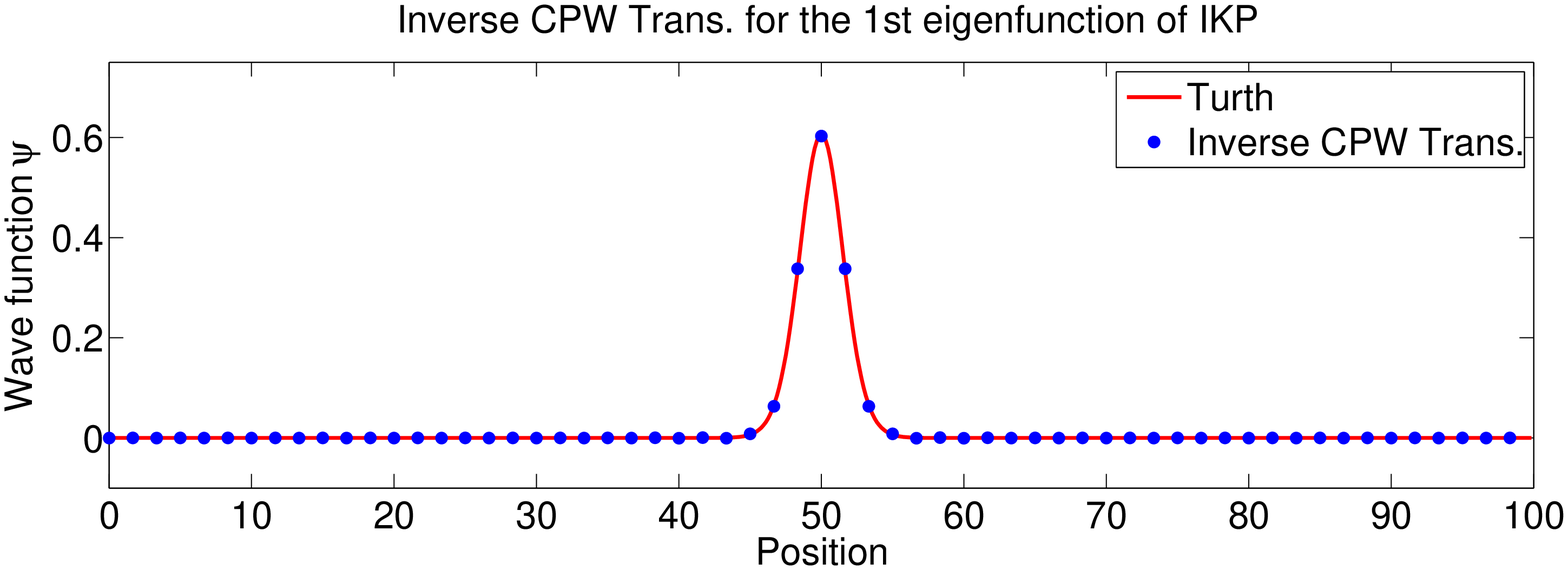}\\
\end{minipage}\hfill
\begin{minipage}{0.49\linewidth}
\centering
\includegraphics[width=1\linewidth]{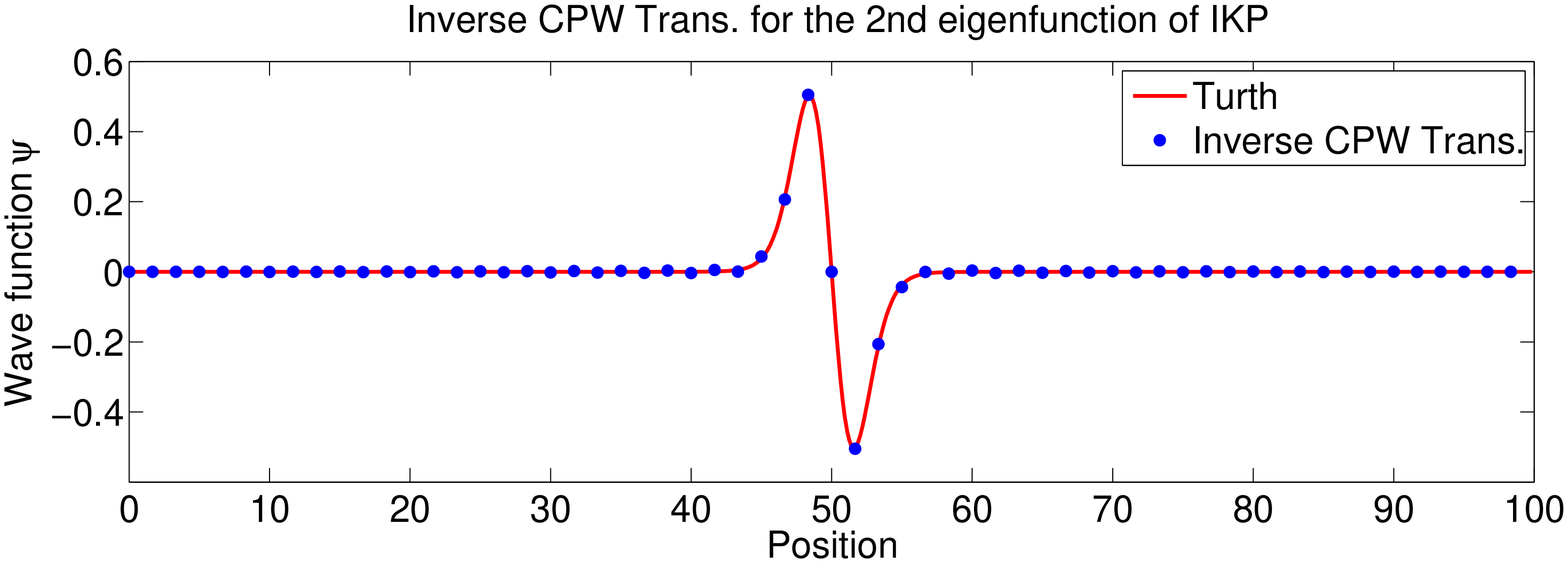}\\
\end{minipage}\hfill\\
\begin{minipage}{0.49\linewidth}
\centering
\includegraphics[width=1\linewidth]{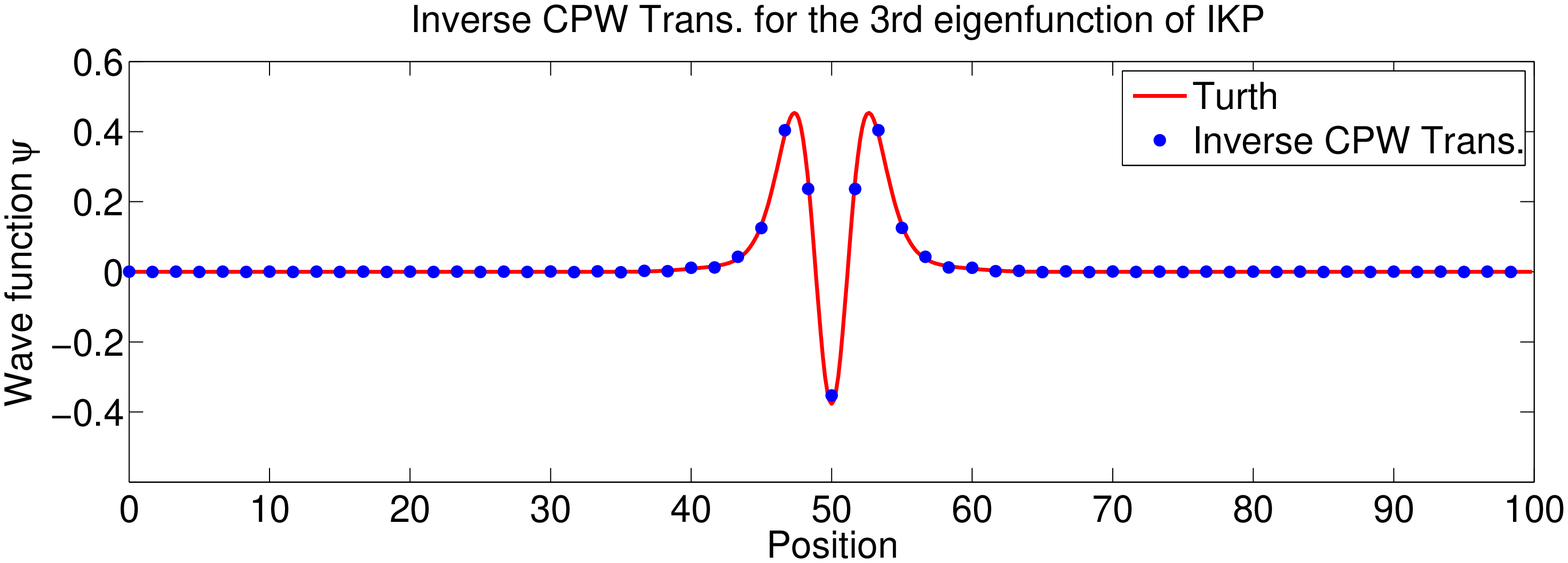}\\
\end{minipage}\hfill
\begin{minipage}{0.49\linewidth}
\centering
\includegraphics[width=1\linewidth]{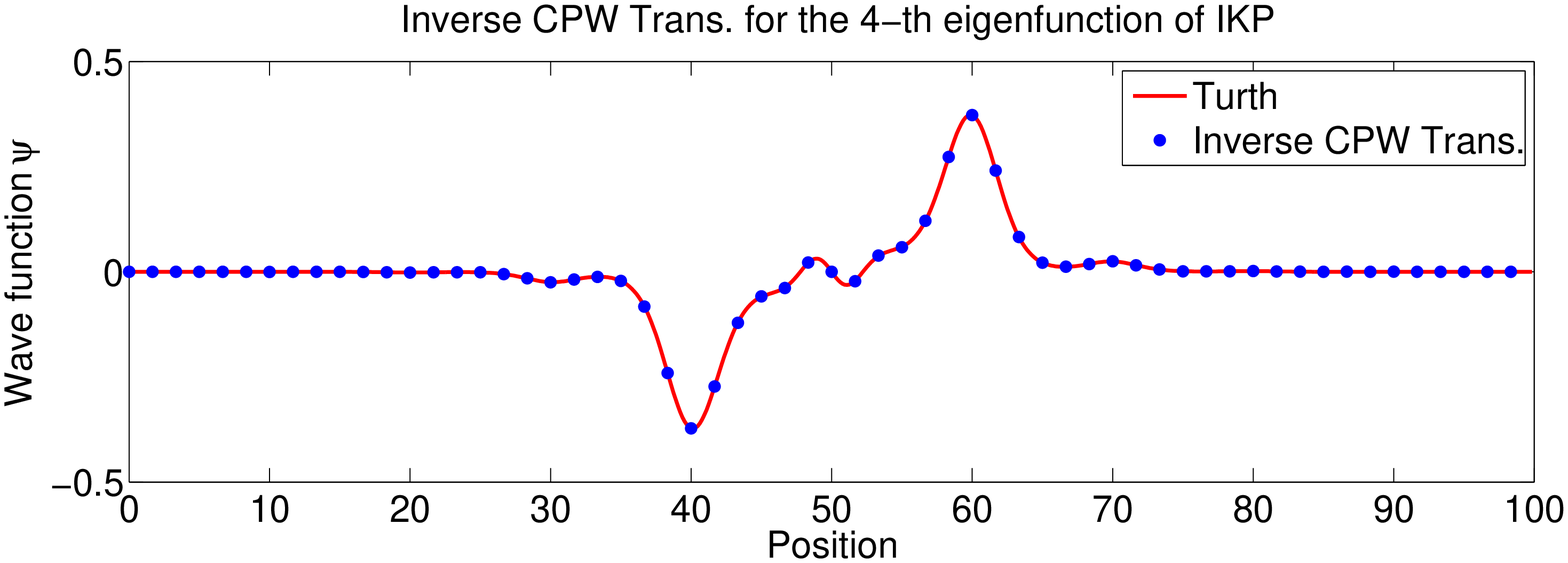}\\
\end{minipage}\hfill
\caption{Reconstruction results using inverse CPW transform for the first four lowest-energy states of the IKP model.}
\label{fig:IKP_InverseCPWtransform}
\end{figure}

Using the first four lowest-energy states of the IKP model illustrated in Figure~\ref{fig:CPW_rep} as examples, we test the inverse CPW transform based on the above algorithm. Figure~\ref{fig:IKP_InverseCPWtransform} shows the results of inverse CPW transforms for the first four energy states of the IKP model, where the solid red line shows the ``exact" results with a very fine mesh, while blue dots are obtained using the fast transform described above with a coarse mesh. 
%Notice that a different number of points is required at each level to get quantitatively accurate representation of the function. Of course, this is consistent with the spectral properties of CPWs, higher $n$ values correspond to higher Fourier components.

\begin{figure}[h]
\begin{minipage}{0.49\linewidth}
\centering
\includegraphics[width=1\linewidth]{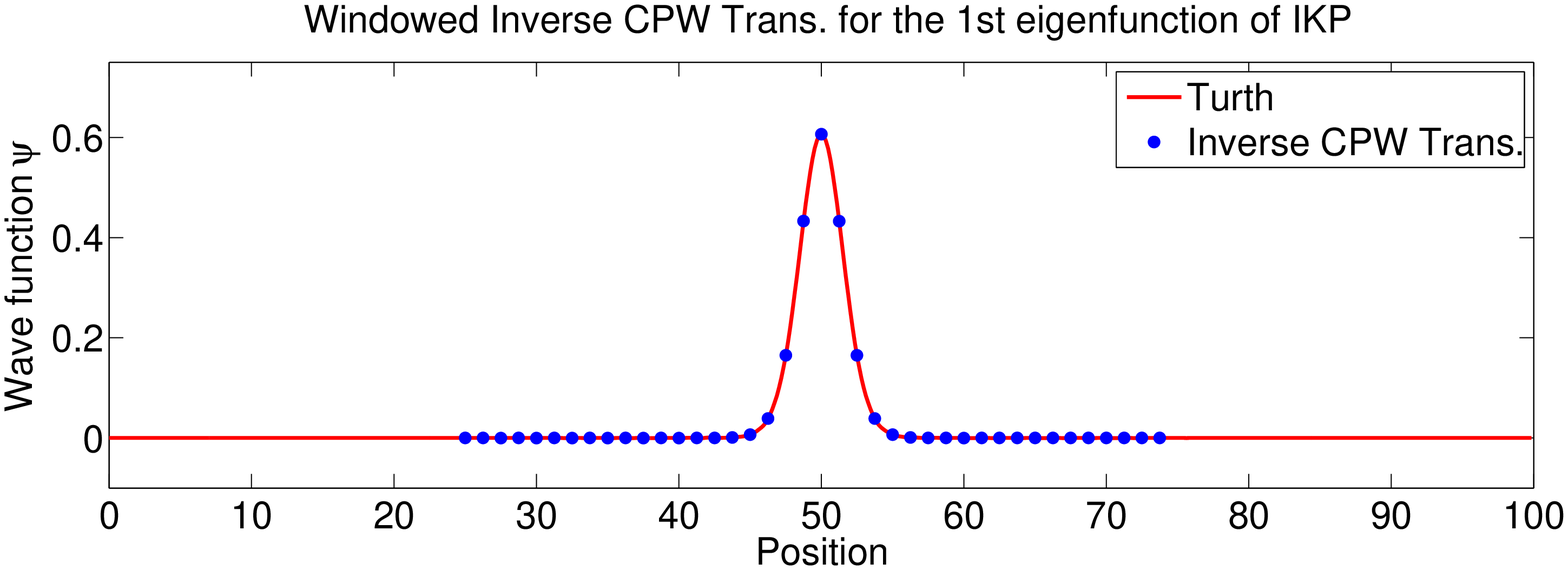}\\
\end{minipage}\hfill
\begin{minipage}{0.49\linewidth}
\centering
\includegraphics[width=1\linewidth]{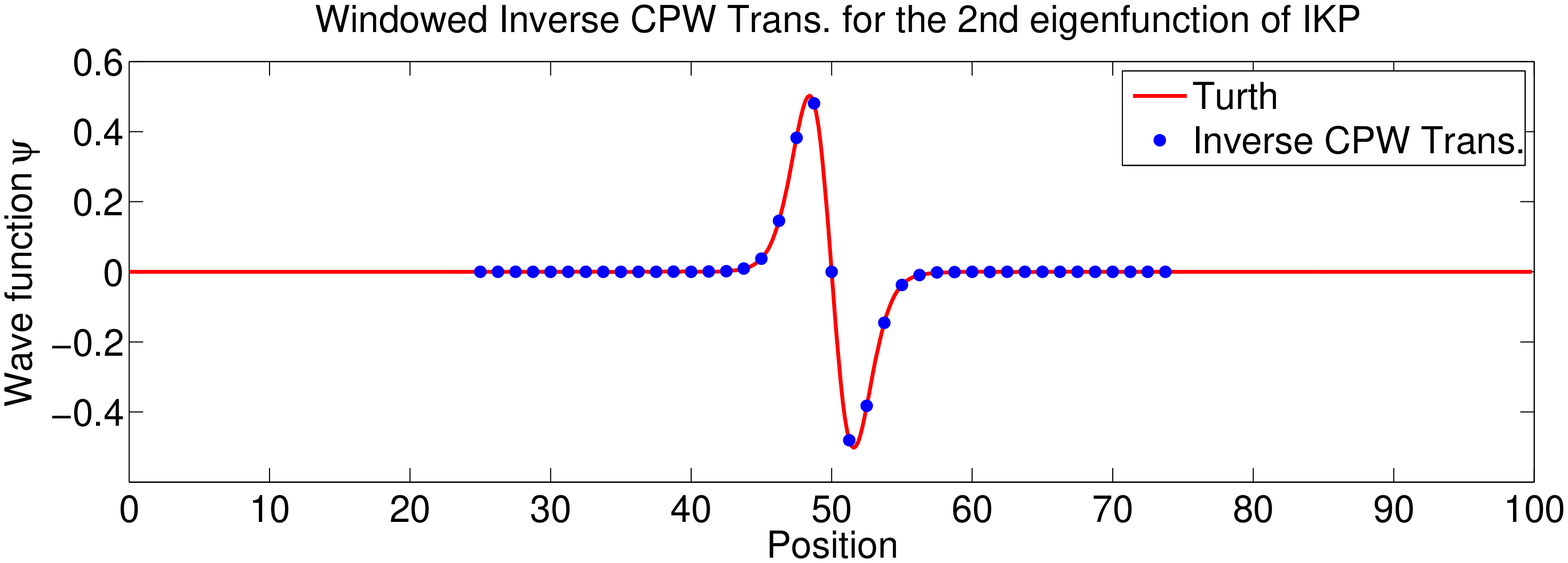}\\
\end{minipage}\hfill\\
\begin{minipage}{0.49\linewidth}
\centering
\includegraphics[width=1\linewidth]{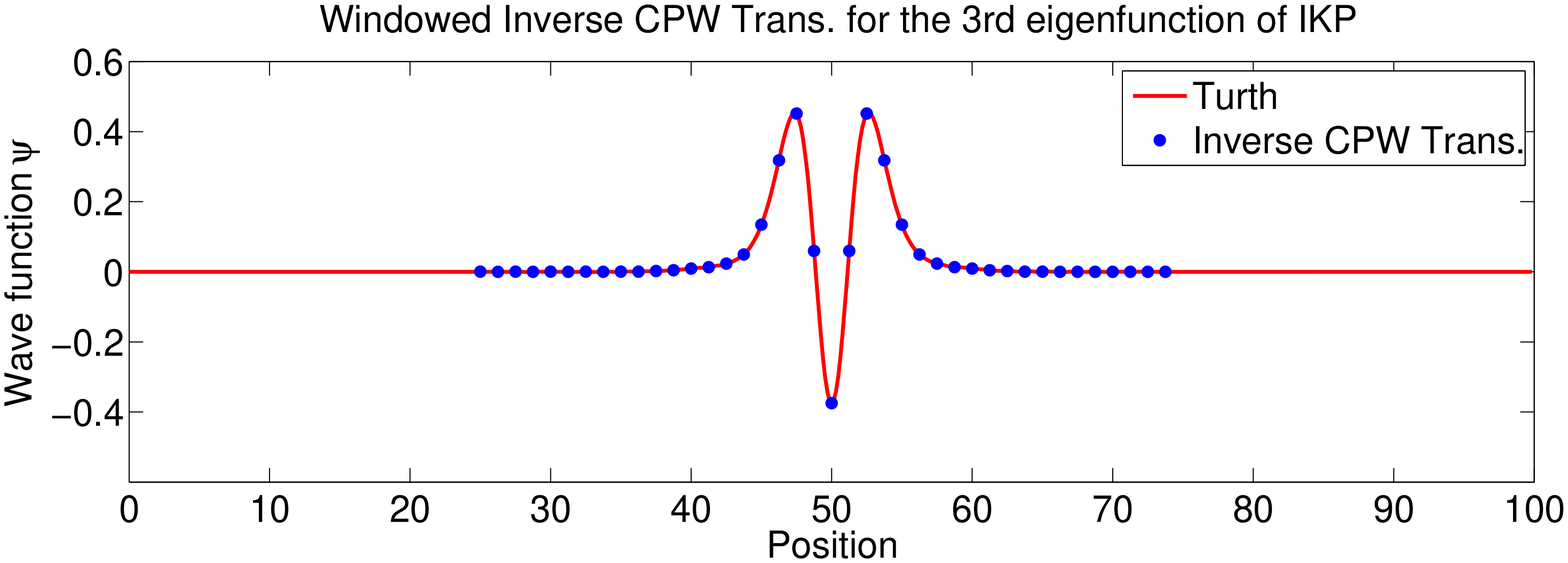}\\
\end{minipage}\hfill
\begin{minipage}{0.49\linewidth}
\centering
\includegraphics[width=1\linewidth]{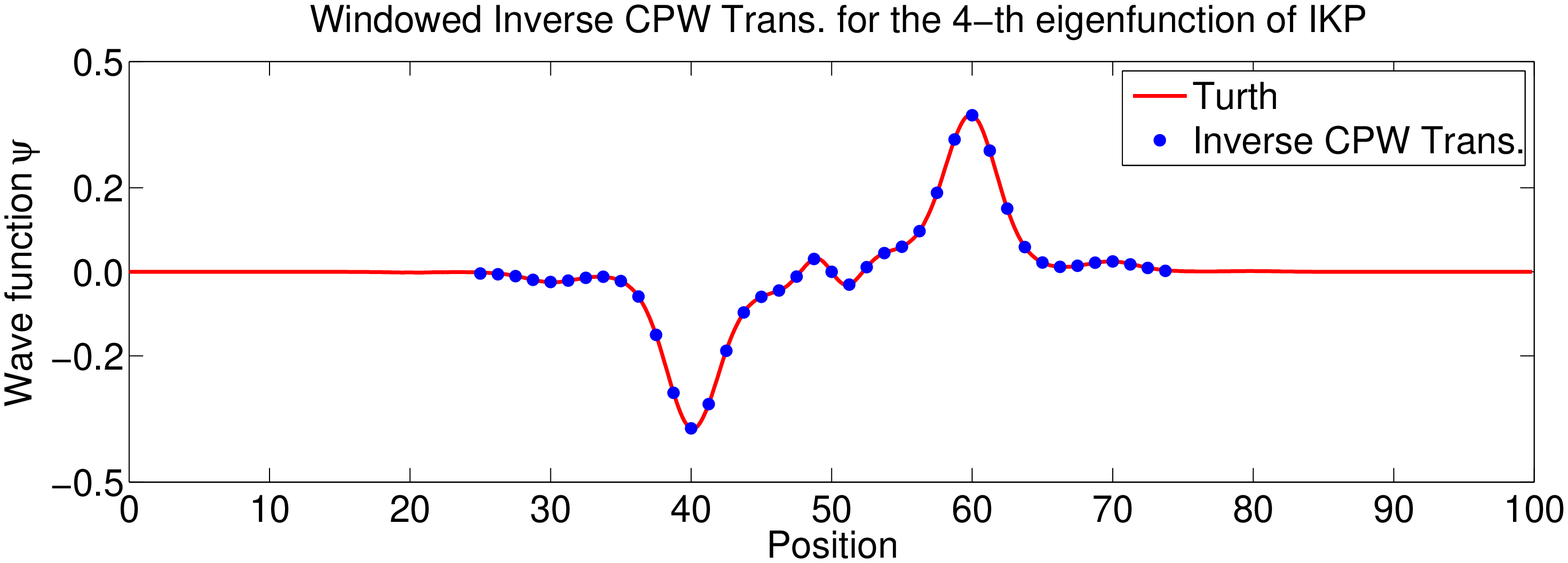}\\
\end{minipage}\hfill
\caption{Reconstruction results using windowed inverse CPW transform for the first four lowest-energy states of the IKP model.}
\label{fig:IKP_ICPWtransform_windowed}
\end{figure}

Finally, these transforms can be ``windowed", which would allow using different number of mesh points in different spatial regions, instead of having to use the same real space mesh everywhere. We conduct numerical tests for the same examples used in Figure~\ref{fig:IKP_InverseCPWtransform} except that these transforms are carried out only over the region where functions are nonzero. Figure~\ref{fig:IKP_ICPWtransform_windowed} reports the results obtained from the ``windowed'' inverse CPW transform, where the solid red line shows the ``exact" results, while the blue dots are obtained using the ``windowed" inverse CPW transform. 
Similarly, the CPW transform can also be windowed just like the ``windowed" inverse CPW transform. 
This ``windowed'' transform will be useful when one needs to use a higher resolution only in certain limited regions.

%%%-------------------------------------------------------------------------------------------------------------------
\section{Discussions and Conclusions}

We have presented a method for producing compressed modes (i.e., modes that are sparse and localized) for the Laplace operator plus a potential $V$, using a variational principle with an $L_1$ penalization term that promotes sparsity.   The SOC algorithm \cite{Lai:2013JSC} was used to numerically construct these modes. In addition, setting $V=0$ we have produced Compressed Plane Waves (CPWs) that form a natural set of modes with two parameters representing position and scale, as in wavelets.
Unlike wavelets, CPWs come from a differential equation so that they have a natural extension to higher dimensions and they may be used as a natural basis for solving PDEs.

These results are only the beginning. We expect to extend CMs techniques in a number of ways and exploit them for a variety of applications:
\begin{enumerate}
  \item Construct CMs for a variety of potentials.
  \item Use the CMs to construct an accelerated (i.e.,  $O(N)$) simulation method for DFT~\cite{Bracewell:1986fourier}.
  \item Develop CMs as the modes for a Galerkin method for PDEs, such as Maxwell's equations.
  \item Generalize CMs for use in PDEs (such as heat type equations) that come from the gradient descent of a variational principle.
  \item Extend CMs to higher dimensions and different geometries, including the Laplace-Beltrami equation on a manifold and a discrete Laplacian on a network. 
\end{enumerate}  

Finally, we plan to perform an analysis of CMs and CPWs to rigorously analyze their existence and qualitative properties, including those that are hypothesized  above.

\section*{Acknowledgments}
V.O. gratefully acknowledges financial support from the 
National Science Foundation under Award Number DMR-1106024 and 
use of computing resources at the National Energy Research Scientific 
Computing Center, which is supported by the US DOE under 
Contract No. DE-AC02-05CH11231. The research of R.C. is partially supported by the US DOE under 
Contract No. DE-FG02-05ER25710. The research of S.O. was supported by the
Office of Naval Research (Grant N00014-11-1-719).

%\bibliographystyle{unsrt}
%% argument is your BibTeX string definitions and bibliography database(s)
%\bibliography{CompressedModes}

\end{document}